\documentclass{amsart}

\usepackage{amsfonts}
\usepackage{amssymb}
\usepackage{latexsym}

\newtheorem{theorem}{Theorem}[section]
\newtheorem{lemma}[theorem]{Lemma}
\newtheorem{definition}[theorem]{Definition}
\newtheorem{proposition}[theorem]{Proposition}

\newtheorem{corollary}[theorem]{Corollary}

\theoremstyle{definition}
\newtheorem{remark}[theorem]{Remark}
\newtheorem{example}[theorem]{Example}

\numberwithin{equation}{section}

\newcommand{\C}{\mathbb{C}}
\newcommand{\Q}{\mathbb{Q}}
\newcommand{\R}{\mathbb{R}}
\newcommand{\Z}{\mathbb{Z}}
\newcommand{\F}{\mathbb{F}}
\newcommand{\Imt}{\mbox{Im}\,}

\begin{document}

\title[General orbifold Euler characteristics and wreath products]
{Generalized orbifold Euler characteristics
for general orbifolds and wreath products}

\author{Carla Farsi}
\address{Department of Mathematics, University of Colorado at Boulder, Campus
Box 395, Boulder, CO 80309-0395 } \email{farsi@euclid.colorado.edu}

\author{Christopher Seaton}
\thanks{The second author was partially supported by a Rhodes College Faculty Development Endowment Grant.}
\address{Department of Mathematics and Computer Science,
Rhodes College, 2000 N. Parkway, Memphis, TN 38112}
\email{seatonc@rhodes.edu}



\begin{abstract}

We introduce the $\Gamma$-Euler-Satake characteristics of a general
orbifold $Q$ presented by an orbifold groupoid $\mathcal{G}$,
generalizing to orbifolds that are not necessarily global quotients
the generalized orbifold Euler characteristics of Bryan-Fulman and
Tamanoi.  Each of these Euler characteristics is defined as the
Euler-Satake characteristic of the space of $\Gamma$-sectors of the
orbifold where $\Gamma$ is a finitely generated discrete group. We
study the behavior of these characteristics under product operations
applied to the group $\Gamma$ as well as the orbifold and establish
their relationships to existing Euler characteristics for orbifolds.
As applications, we generalize formulas of Tamanoi, Wang, and Zhou
for the Euler characteristics and Hodge numbers of wreath symmetric
products of global quotient orbifolds to the case of quotients by
compact, connected Lie groups acting almost freely.

\end{abstract}

\maketitle


\section{Introduction}
\label{sec-intro}

When Satake first introduced orbifolds under the name
\emph{$V$-manifolds}, one of the first invariants defined was the
Euler-Satake characteristic, then called the \emph{Euler
characteristic as a $V$-manifold}; see \cite{satake2}.  Since that
time, a number of Euler characteristics have been introduced for
orbifolds.  Most notably, the ``stringy" orbifold Euler
characteristic was introduced in \cite{dixon} for global quotients
and later generalized to general orbifolds in \cite{roan}; see also
\cite{atiyahsegal} and \cite{hirzebruchhoefer}.  These Euler
characteristics were generalized for global quotients twice---first
by Bryan and Fulman in \cite{bryanfulman}, where they appeared as
the first two elements of a sequence of Euler characteristics, and
independently by Tamanoi in \cite{tamanoi1} and \cite{tamanoi2},
where an Euler characteristic was defined for each group $\Gamma$,
the sequence of Bryan and Fulman corresponding to $\Gamma = \Z^m$.

Here, we extend these definitions, introducing the
\emph{$\Gamma$-Euler-Satake characteristics} for a general orbifold
$Q$ presented by an orbifold groupoid $\mathcal{G}$.  Given a
finitely generated discrete group $\Gamma$, we define the
$\Gamma$-Euler-Satake characteristic as the Euler-Satake
characteristic of the \emph{$\Gamma$-sectors of $Q$}, an orbifold generalizing
the inertia orbifold and multi-sectors (see \cite{farsiseaton1} and \cite{farsiseaton2}).

Our main results are as follows. Generalizing a result of
\cite{tamanoi1}, Theorem \ref{thrm-Gamma1timesGamma2formula} shows
that recursively constructing $\Gamma$-sectors corresponds to the direct
product operation for groups.  This allows us to relate the
$\Gamma$-Euler-Satake characteristics to the usual (topological) Euler
characteristic of the underlying space of the orbifold in
question. We then demonstrate Theorem \ref{thrm-mainWPformula},
establishing a formula for the $\Gamma$-Euler-Satake characteristic
of a wreath symmetric product by a compact, connected Lie group that
generalizes results of \cite{wang1}.  Theorem
\ref{thrm-hodgeproductformula} generalizes a formula for the
(shifted) Hodge numbers of wreath symmetric products from the case
of finite group actions to the case of compact, complex, connected
Lie group actions.

In Section \ref{sec-prelim}, we recall the relevant preliminary
material and fix notation.  In particular, in Subsection
\ref{subsec-eulersatakecharproperties}, we verify the
multiplicativity and other basic properties of the Euler-Satake
characteristic, and in Subsections \ref{subsec-wreathprodcompact}
and \ref{subsec-actionswreathprodmanifolds}, we discuss the
structures of wreath products of compact, connected Lie groups and
the associated wreath symmetric product orbifolds.  Section
\ref{sec-operationsgroupsorbifolds} establishes the behavior the
$\Gamma$-sectors of a general orbifold under product operations on
both the orbifold and the group $\Gamma$.  In Section
\ref{sec-GammaESC}, we define the $\Gamma$-Euler-Satake
characteristics and demonstrate their connections with other Euler
characteristics for orbifolds.  In Section \ref{sec-wreathproducts},
we turn our attention to wreath symmetric products by compact, connected Lie groups and prove Theorem \ref{thrm-mainWPformula}.
Section \ref{sec-hodgenumbers} turns to the case of complex
orbifolds given by quotients by compact, connected groups, proving
Theorem \ref{thrm-hodgeproductformula}.  Note that throughout,
we assume that $G$ acts effectively for simplicity; using
Equation \ref{eq-ESCnoneff} below, our results can be applied the noneffective case.

This paper presents the first of our investigation of these
generalized Euler characteristics for orbifolds.  In future work, we
will continue these investigations, expecting to extend results in
\cite{tamanoi2}, \cite{tamanoi3}, \cite{wang1}, and \cite{wang2} to
more general classes of orbifolds and further explore group actions
on orbifolds.

The first author would like to thank the MSRI for its hospitality
during the preparation of this manuscript.


\section{Preliminaries and Background Material}
\label{sec-prelim}

In this section, we briefly recall the definitions we will need. For
background on orbifolds, the reader is referred to
\cite{ademleidaruan}; see also \cite{moerdijkmrcun},
\cite{moerdijkorbgroupintro}, \cite{chenruangwt} or \cite{ruansgt}.
For details on the definition of the $\Gamma$-sectors of an
orbifold, see \cite{farsiseaton1} and \cite{farsiseaton2}.


\subsection{Orbifolds and $\Gamma$-Sectors}
\label{subsec-orbifoldsectordefs}

An \emph{orbifold structure} on a paracompact Hausdorff space $\mathbb{X}_Q$ is an orbifold groupoid
$\mathcal{G}$, i.e. a proper \'{e}tale Lie groupoid, and a homeomorphism
$f: |\mathcal{G}| \rightarrow \mathbb{X}_Q$ between the orbit space
$|\mathcal{G}|$ of $\mathcal{G}$ and $\mathbb{X}_Q$.  We say that $(\mathcal{G}, f)$ is a \emph{presentation} of the orbifold structure.  Two presentations $(\mathcal{G}, f)$ and $(\mathcal{G}^\prime, f^\prime)$
are \emph{equivalent} if $\mathcal{G}$ and $\mathcal{G}^\prime$ are Morita equivalent and the homeomorphisms
$g : |\mathcal{G}| \rightarrow |\mathcal{G}^\prime|$ induced by the Morita equivalence satisfies
$f = f^\prime \circ g$.  An \emph{orbifold} $Q$ is a paracompact Hausdorff space $\mathbb{X}_Q$, called
the \emph{underlying space of $Q$}, and an equivalence class of orbifold structures on $\mathbb{X}_Q$.
Given a presentation $(\mathcal{G}, f)$,
of the orbifold $Q$, we will often identify $\mathbb{X}_Q$ with $|\mathcal{G}|$ and avoid explicit reference to $f$.  We say that orbifolds $Q_1$ and $Q_2$ are
\emph{orbifold-diffeomorphic} or simply \emph{diffeomorphic} if the
groupoids representing them are Morita equivalent, so that in particular their underlying spaces are homeomorphic.

Let $Q$ be an orbifold.  Throughout, we use the notation that
$\mathcal{G}$ is an orbifold groupoid presenting $Q$ with space of
objects $G_0$ and space of arrows $G_1$.  When considering wreath
symmetric products, we will restrict our attention to orbifolds $Q$
presented by $M \rtimes G$ where $G$ is a compact, connected Lie
group acting smoothly, effectively, and locally freely (i.e.
properly with discrete stabilizers) on the smooth manifold $M$ so
that $M\rtimes G$ is Morita equivalent to an orbifold groupoid. An
orbifold presented in such a way is called a \emph{quotient
orbifold}, and a \emph{global quotient orbifold} if $G$ is finite.
Note that $M\rtimes G$ is not \'{e}tale unless $G$ is finite.  For
clarity, we will represent general orbifolds in terms of left
groupoid actions and quotient orbifolds in terms of right group
actions.

For every $x \in G_0$, there is an open neighborhood $V_x \subseteq G_0$ of $x$ diffeomorphic to $\R^n$
with $x$ corresponding to the origin such that the isotropy group $G_x$ acts linearly
on $V_x$ and the restriction $\mathcal{G}|_{V_x}$ is isomorphic to $G_x \ltimes V_x$.
We let $\pi_x : V_x \rightarrow |G_x \ltimes V_x| \subseteq |\mathcal{G}|$ denote
the quotient map, and refer to the triple $\{ V_x, G_x, \pi_x \}$ as a \emph{linear orbifold chart}
for $Q$ at $x$.

Let $\mathcal{G}$ and $\mathcal{H}$ be orbifold groupoids.  By
$\mathcal{G} \times \mathcal{H}$, we mean the groupoid with objects
$G_0 \times H_0$ and arrows $G_1 \times H_1$ (see \cite[page
123]{moerdijkmrcun}).  As the products of proper maps is proper and
the product of local diffeomorphism is a local diffeomorphism,
$\mathcal{G} \times \mathcal{H}$ is an orbifold groupoid.  Linear charts
for this orbifold are given by products of linear charts for
$|\mathcal{G}|$ and $|\mathcal{H}|$.

Given a finitely generated discrete group $\Gamma$, the space
$\mathcal{S}_\mathcal{G}^\Gamma = \mbox{HOM}(\Gamma, \mathcal{G}) = \bigcup_{x \in G_0} \mbox{HOM}(\Gamma, G_x)$ of
groupoid homomorphisms from $\Gamma$ into $\mathcal{G}$ naturally inherits the structure of a smooth
$\mathcal{G}$-manifold.  We let $\mathcal{G}^\Gamma = \mathcal{G}\ltimes \mathcal{S}_\mathcal{G}^\Gamma$
denote the corresponding translation groupoid, an orbifold groupoid, and $\tilde{Q}_\Gamma$ the
corresponding orbifold.  We refer to $\tilde{Q}_\Gamma$ as the \emph{orbifold of $\Gamma$-sectors
of $Q$}.  If $\phi_x : \Gamma \rightarrow G_x$ is an element of $\mathcal{S}_\mathcal{G}^\Gamma$, we let
$\tilde{Q}_{(\phi)}$ denote connected component of $\tilde{Q}_\Gamma$ containing the orbit of $\phi_x$ and refer to $\tilde{Q}_{(\phi)}$ as the \emph{$\Gamma$-sector associated to $\phi_x$}.

Given a homomorphism $\phi_x : \Gamma \rightarrow G_x$ and
a linear chart $\{ V_x, G_x, \pi_x \}$ for $Q$ near $x$, there is a diffeomorphism
$\kappa_{\phi_x} : V_x^{\langle \phi_x \rangle} \rightarrow \mathcal{S}_\mathcal{G}^\Gamma$ onto
a neighborhood of $\phi_x$, where $V_x^{\langle \phi_x \rangle}$ denotes the subspace of $V_x$
fixed by the image of $\phi_x$.  Then up to identification via $\kappa_{\phi_x}$,
$\{ V_x^{\langle \phi_x \rangle}, C_{G_x}(\phi_x), \pi_x^{\phi_x} \}$ forms a linear chart for
the groupoid $\mathcal{G}^\Gamma = \mathcal{G}\ltimes \mathcal{S}_\mathcal{G}^\Gamma$ at $\phi_x$
where $C_{G_x}(\phi_x)$ is the centralizer of $\Imt \phi_x$ in $G_x$ and $\pi_x^{\phi_x}$ is the
quotient map.

In the case that $Q$ is presented by $M \rtimes G$ where $M$ is a
smooth manifold and $G$ is a compact Lie group acting smoothly and
locally freely, the space of $\Gamma$-sectors admits the
presentation
\[
    (M;G)_\Gamma := \coprod\limits_{(\phi) \in t_{M;G}^\Gamma}
    M^{\langle \phi \rangle} \rtimes C_G(\phi)
\]
where $t_{M;G}$ denotes the set of conjugacy classes of homomorphisms
$\phi \in \mbox{HOM}(\Gamma, G)$ such that $M^{\langle \phi \rangle} \neq \emptyset$.
Note that $M^{\langle \phi \rangle} \rtimes C_G(\phi)$ need not be a connected orbifold;
see \cite[Definition 2.1, Theorem 3.6]{farsiseaton2} for details.

\subsection{Cohomology and $K$-Theory}
\label{subsec-cohomology}

Unless otherwise specified, all cohomology is with complex coefficients.

If $Q$ is an orbifold presented by the orbifold groupoid
$\mathcal{G}$, we let $H^\ast(\mathcal{G})$ or $H^\ast(Q)$ denote
the (complex) de Rham cohomology of $\mathcal{G}$-invariant
differential forms on $G_0$.  Note that $H^\ast(\mathcal{G})$ is
isomorphic to the usual singular cohomology of the underlying space
of $Q$ as well as the cohomology $H_{orb}^\ast(Q)$ defined in terms
of groupoid classifying spaces; see \cite[Section
2.1]{ademleidaruan}. Hence, we need not distinguish between
$H^\ast(\mathcal{G})$ and $H^\ast(Q)$.

If $\Gamma$ is a finitely generated discrete group, then we let
$H_\Gamma^\ast(Q) = H^\ast\left( \tilde{Q}_\Gamma \right)$ denote
the de Rham cohomology of the orbifold of $\Gamma$-sectors
$\tilde{Q}_\Gamma$, called the $\Gamma$-cohomology.  In particular,
it follows from \cite[Corollary 3.8]{farsiseaton2} that the
$\Z$-sectors of $Q$ are orbifold-diffeomorphic to the inertia
orbifold of $Q$. Hence, $H_\Z^\ast(Q) = H^\ast\left( \tilde{Q}_\Z
\right)$ is the \emph{delocalized cohomology of $Q$}.  As additive
groups, the delocalized cohomology corresponds to the Chen-Ruan
orbifold cohomology; see \cite{ademleidaruan} or
\cite{chenruanorbcohom}.

If $Q$ is a quotient orbifold, then the orbifold $K$-theory of $Q$
can be defined in a variety of ways; see \cite{baumconnes},
\cite{farsi}, \cite{tuxu}, or \cite{ademruan}.  When complexified,
each of these definitions is isomorphic, and is isomorphic via the
Chern character to the delocalized cohomology of $Q$.

In the case that $Q$ is a complex orbifold, i.e. an orbifold
presented by a holomorphic orbifold groupoid $\mathcal{G}$ where
$G_0$ is a complex manifold and the elements of $G_1$ are locally
automorphisms, we define the Dolbeault cohomology
$H^{\ast,\ast}(\mathcal{G}) = H^{\ast,\ast}(Q)$ of $Q$ in terms of
$G_1$-invariant forms on $G_0$; see \cite{ademleidaruan}.

It follows from the same argument as \cite[Corollary
4.2]{ademleidaruan} for multi-sectors that each $\Gamma$-sector
inherits a complex structure.  Hence $\tilde{Q}_\Gamma$ is a complex
orbifold, and we define the $\Gamma$-Dolbeault cohomology as
\[
    H_\Gamma^{\ast,\ast}(Q)
    =
    H_\Gamma^{\ast,\ast} \left( \tilde{Q}_\Gamma \right).
\]
In particular, if $G$ is a compact, complex Lie group that acts
effectively, locally freely, and holomorphically on the compact
complex manifold $M$, then $M \rtimes G$ is a complex orbifold.  As
the action is locally free, the fixed-point set of each subgroup $H
\leq G$ is either empty or the intersection of a finite collection
of complex submanifolds.  It follows that for each $H \leq G$ such
that $M^H \neq \emptyset$, $M^H$ is a complex submanifold, and the
restriction of the $C_G(H)$-action to $M^H$ is holomorphic.


\subsection{The Euler-Satake Characteristic and its Properties}
\label{subsec-eulersatakecharproperties}

Let $Q$ be a compact orbifold of dimension $n$ and let
$Q_{\mbox{\scriptsize eff}}$ be the effective orbifold associated to
$Q$ (see \cite[Definition 2.33]{ademleidaruan}).  It is well known
that $Q_{\mbox{\scriptsize eff}}$ can be presented by a groupoid $M
\rtimes G$ where $M$ is a manifold and $G$ a compact Lie group
acting smoothly and locally freely on $M$. By \cite[page
488]{illman} (see also \cite{yang} and \cite{verona}), it follows
that $Q_{\mbox{\scriptsize eff}}$ admits a good finite
triangulation; i.e. a finite triangulation in which the $G$-isotropy
type is constant on the interior of each simplex.  Note that if
$c:[0,1]\rightarrow Q$ is a curve in $Q$, then the isomorphism type
of the isotropy group for points in $c$ is constant if and only if
it is constant on the induced curve in $Q_{\mbox{\scriptsize eff}}$.
It follows that a good triangulation of $Q_{\mbox{\scriptsize eff}}$
induces a triangulation of $Q$ such that the isomorphism type of the
isotropy group is constant on the interior of each simplex. By
refining triangulations and using stellar subdivisions if necessary,
we may also assume that each of the top simplices is contained
inside the image of a linear orbifold chart; see
\cite{moerdijkpronksimplicial}. Following the language for
quotients, we will refer to such a triangulation of the orbifold $Q$
as \emph{good}.

The following definition was originally stated in \cite{satake2} under the name \emph{Euler characteristic as a $V$-manifold}.

\begin{definition}[The Euler-Satake Characteristic]
\label{def-ESC}

Let $Q$ be a closed orbifold and $\mathcal{T}$ a good triangulation of $Q$.  The \emph{Euler-Satake characteristic of $Q$} is
\[
    \chi_{ES}(Q)
    =
    \sum\limits_{\sigma \in \mathcal{T}} (-1)^{\mbox{\scriptsize dim}\: \sigma}
    \frac1{|G_{\sigma}|},
\]
where $G_\sigma$ denotes the isotropy group of a point on the
interior of $\sigma$. If $Q^\prime$ is a subset of $Q$ corresponding
to a subcomplex $\mathcal{T}^\prime$ of $\mathcal{T}$, then we
define $\chi_{ES}(Q^\prime)$ identically, summing over those
simplices contained in $\mathcal{T}^\prime$.

\end{definition}

Note that if $\mathcal{G}$ is a groupoid presenting $Q$, it will
frequently be convenient for us to use the notation
$\chi_{ES}(\mathcal{G})$ to denote $\chi_{ES}(Q)$.

\subsubsection{Properties of the Euler-Satake Characteristic}
\label{subsubsec-ESCproperties}

It is clear that $\chi_{ES}(Q)$ does not depend on the choice of
(good) triangulation. By contrast, the standard Euler characteristic
$\chi_{top}(Q)$ of the underlying topological space of $Q$ can be
expressed as
\[
    \chi_{top}(Q) =   \sum\limits_{\sigma \in \mathcal{T}} (-1)^{\mbox{\scriptsize dim}\: \sigma}.
\]
Hence, when $Q$ is a manifold, i.e. when $G_\sigma =1$ for each
$\sigma \in \mathcal{T}$, $\chi_{ES}(Q) = \chi_{top}(Q)$.

Note that if $Q$ is a global quotient orbifold presented by $M\rtimes G$ so that $G$ is a finite group, then
\begin{equation}
\label{eq-ESCglobquot}
    \chi_{ES}(Q)
    =
    \frac{1}{|G|}\chi_{top}(M).
\end{equation}
Similarly, it is easy to see that if $Q$ is connected and noneffective, then
\begin{equation}
\label{eq-ESCnoneff}
    \chi_{ES}(Q)
    =
    \frac{1}{|K_p|}\chi_{ES}(Q_{\mbox{\scriptsize eff}})
\end{equation}
where $K_p$ is the isotropy group of any point $p$ that is
nonsingular in $Q_{\mbox{\scriptsize eff}}$ (equivalently the normal
subgroup $K_p \leq G_p$ of any isotropy group that acts trivially in
a chart).  In the case that $Q$ is not connected, the isomorphism
class of $K_p$ may vary over connected components.

It is also a direct consequence of the definition that the Euler-Satake characteristic is additive;
that is, if $Q_1$ and $Q_2$ are subsets of $Q$ such that $Q_1 \cup Q_2 = Q$, and the sets $Q_1$,
$Q_2$, and $Q_1 \cap Q_2$ correspond to subcomplexes of $\mathcal{T}$, then
\begin{equation}
\label{eq-ESCadditive}
    \chi_{ES}(Q_1 \cup  Q_2)
    =
    \chi_{ES}(Q_1) + \chi_{ES}(Q_2)
    - \, \chi_{ES}(Q_1 \cap Q_2).
\end{equation}

A \emph{covering orbifold} of the orbifold $Q$ presented by the
orbifold groupoid $\mathcal{G}$ is a $\mathcal{G}$-space $E$
equipped with a connected covering projection $\rho: E \rightarrow
G_0$; see \cite[page 40]{ademleidaruan}.  Note that this definition
generalizes \cite[Definition 13.2.2]{thurston} to the case of
noneffective orbifolds, and further that it requires that the
kernels of the actions on $Q$ and the cover coincide.  It follows
that if $\rho: \hat{Q} \rightarrow Q$ is a covering orbifold, then
$\rho_{\mbox{\scriptsize eff}} : \hat{Q}_{\mbox{\scriptsize eff}}
\rightarrow Q_{\mbox{\scriptsize eff}}$ is a covering orbifold as
well. We define the \emph{number of sheets} of the covering $\rho :
\hat{Q} \rightarrow Q$ to be the number of points in a preimage of a
nonsingular point of $Q_{\mbox{\scriptsize eff}}$. We have the
following.

\begin{lemma}
\label{lem-ESCmultipliconcovers}

Let $\rho: \hat{Q} \rightarrow Q$ be a covering orbifold with $k$
sheets.  Then
\[
    \chi_{ES}(\hat{Q})
    =
    k\chi_{ES}(Q).
\]

\end{lemma}

See \cite[Proposition 13.3.4]{thurston} for the effective case, and
note that the noneffective case follows from Equation
\ref{eq-ESCnoneff}.

\subsubsection{Multiplicativity of the Euler-Satake Characteristic}
\label{subsubsec-ESCmultiplic}

If $\mathcal{T}_1$ and $\mathcal{T}_2$ are triangulations of spaces
$Q_1$ and $Q_2$, respectively, then the simplicial product
$\mathcal{T}_1 \times \mathcal{T}_2$ is a triangulation of $Q_1
\times Q_2$; see \cite[page 61]{nfryv} for the construction.  We
have the following.

\begin{proposition}
\label{prop-productgoodsimplicialcomplexgood}

The simplicial product of two good triangulations is good.

\end{proposition}

\begin{proof}

It is enough to notice that, by construction, the interior of a
simplex in the simplicial product is contained in the product of
interiors of simplices in $\mathcal{T}_1$ and $\mathcal{T}_2$. It
follows that if $p \in \sigma_1 \in \mathcal{T}_1$ and $q \in
\sigma_2 \in \mathcal{T}_2$, then $(p, q) \in Q_1 \times Q_2$ has
isotropy $G_p \oplus G_q$.  Moreover, if $\sigma_1 \in
\mathcal{T}_1$ is contained in the image of the linear chart $\{
V_1, G_1, \pi_1 \}$ for $Q_1$ and $\sigma_2 \in \mathcal{T}_2$ is
contained in the image of the linear chart $\{ V_2, G_2, \pi_2 \}$
for $Q_2$, then $\sigma_1 \times \sigma_2$ is contained in the
linear chart $\{ V_1 \times V_2, G_1 \oplus G_2, \pi_1 \times \pi_2
\}$ for $Q_1 \times Q_2$.

\end{proof}

The main result of this section is the following.

\begin{theorem}[Multiplicativity of the Euler-Satake Characteristic]
\label{thrm-ESCmultiplic}

Let $Q_1$ and $Q_2$ be compact orbifolds.  Then
\[
    \chi_{ES}(Q_1 \times Q_2) = \chi_{ES}(Q_1 ) \chi_{ES}(Q_2).
\]

\end{theorem}

We first prove the following, whose proof loosely follows \cite[page
260]{kawakubo}.

\begin{lemma}
\label{lem-ESCmultiplicmanifold}

Let $Q_1$ be a smooth, compact manifold (possibly with boundary) and $Q_2$ a compact orbifold.  Then
\[
    \chi_{ES}(Q_1 \times Q_2)
    =   \chi_{ES}(Q_1 ) \chi_{ES}(Q_2).
\]

\end{lemma}

Note that $\chi_{ES}(Q_1) = \chi_{top}(Q_1)$.

\begin{proof}

Let $n$ denote the dimension of $Q_2$, let $\mathcal{T}_1$ be a
finite triangulation of $Q_1$, and let $\mathcal{T}_2$ be a good
finite triangulation of $\mathcal{T}_2$.  For $i = 1, 2$, we let
$\mathcal{T}_i^t$ denote the $t$-skeleton of $\mathcal{T}_i$.  Fix
$\sigma^t \in \mathcal{T}_2$, and let $\{ V_{\sigma^t},
G_{\sigma^t}, \pi_{\sigma^t} \}$ be a linear orbifold chart for
$Q_2$ whose image contains $\sigma^t$. Let $\tilde{\sigma}^t =
\pi^{-1}(\sigma^t)$, and then $\pi_{\sigma^t} : \tilde{\sigma}^t
\rightarrow \sigma^t$ is an orbifold cover.  It follows that $Q_1
\times \sigma^t$ is presented by $(Q_1 \times \tilde{\sigma}^t)
\rtimes G_{\sigma^t}$ with $G_{\sigma^t}$ acting trivially on the
first factor.  Hence,
\[
\begin{array}{rcl}
    \chi_{ES}(Q_1 \times \sigma^t)
        &=& \chi_{ES}((Q_1 \times \tilde{\sigma}^t) \rtimes G_{\sigma^t})            \\\\
        &=& \chi_{top}(Q_1 \times \tilde{\sigma}^t)/|G_{\sigma^t}|          \\\\
        &=& \chi_{top}(Q_1)  \chi_{top}(\tilde{\sigma}^t)/|G_{\sigma^t}|          \\\\
        &=& \chi_{ES}(Q_1)  \chi_{ES}(\sigma^t)
\end{array}
\]
by Equation \ref{eq-ESCglobquot}, the multiplicative property of $\chi_{top}$,
and the multiplicative property of $\chi_{ES}$ on orbifold covers.

Applying this to each $\sigma \in \mathcal{T}_2$ and using Equation \ref{eq-ESCadditive}, we have
\[
\begin{array}{rcl}
    \chi_{ES}(Q_1 \times Q_2)
    &=&
    \sum\limits_{\sigma^n \in \mathcal{T}_2^n} \chi_{ES}( Q_1 \times \sigma^n )
    - \sum\limits_{\sigma^{n-1} \in \mathcal{T}_2^{n-1}} \chi_{ES}( Q_1 \times \sigma^{n-1} )
    \\\\ &&
    + \sum\limits_{\sigma^{n-2} \in \mathcal{T}_2^{n-2}} \chi_{ES}( Q_1 \times \sigma^{n-2} )
    + \cdots + (-1)^n
    \sum\limits_{\sigma^0 \in \mathcal{T}_2^0} \chi_{ES}( Q_1 \times \sigma^0 )
            \\\\
    &=&
    \sum\limits_{t=0}^n (-1)^{n - t}
    \sum\limits_{\sigma^t \in \mathcal{T}_2^t} \chi_{ES}( Q_1 \times \sigma^t )
            \\\\
    &=&
    \sum\limits_{t=0}^n (-1)^{n - t}
    \sum\limits_{\sigma^t \in \mathcal{T}_2^t} \chi_{ES}(Q_1)\chi_{ES}(\sigma^t)
            \\\\
    &=&
    \chi_{ES}(Q_1) \sum\limits_{t=0}^n (-1)^{n - t}
    \sum\limits_{\sigma^t \in \mathcal{T}_2^t}\chi_{ES}(\sigma^t)
            \\\\
    &=&
    \chi_{ES}(Q_1)\chi_{ES}(Q_2).
\end{array}
\]

\end{proof}

\begin{proof}[Proof of Theorem \ref{thrm-ESCmultiplic}]

Let $\mathcal{T}_1$ and $\mathcal{T}_2$ be good finite
triangulations of $Q_1$ and $Q_2$, respectively.  We use the same notation as in the proof of
Lemma \ref{lem-ESCmultiplicmanifold} for skeletons and charts.  Applying Lemma
\ref{lem-ESCmultiplicmanifold} to each $\sigma^t \in \mathcal{T}_2$ as well as
multiplicativity on orbifold covers, we have
\[
\begin{array}{rcl}
    \chi_{ES}(Q_1 \times \sigma^t)
    &=&
    \chi_{ES}((Q_1\times \tilde{\sigma}^t)\rtimes G_{\tilde{\sigma}^t})
                \\\\
    &=&
    \chi_{ES}(Q_1 \times \tilde{\sigma}^t)/|G_{\tilde{\sigma}^t}|
                \\\\
    &=&
    \chi_{ES}(Q_1)\chi_{ES}(\tilde{\sigma}^t)/|G_{\tilde{\sigma}^t}|
                \\\\
    &=&
    \chi_{ES}(Q_1)\chi_{ES}(\sigma^t).
\end{array}
\]
Summing up over all of the simplices in $\mathcal{T}_2$ using Equation \ref{eq-ESCadditive}
as in the proof of Lemma \ref{lem-ESCmultiplicmanifold}, it follows that
\[
    \chi_{ES}(Q_1 \times Q_2)
    =
    \chi_{ES}(Q_1)\chi_{ES}(Q_2).
\]

\end{proof}


\subsection{Wreath Products of Compact Groups}
\label{subsec-wreathprodcompact}

In this section, we generalize some of the results of
\cite{tamanoi1} to compact, connected Lie groups acting effectively
and almost freely; see also \cite{wang1}. Throughout this section,
we take $G$ to be a compact, connected Lie group and let $G^n =
\prod_{i=1}^n G$ denote the direct product with elements denoted
$\mathbf{g} = (g_1, \ldots, g_n)$.

Let $\mathcal{S}_n$ denote the symmetric group on a set of $n$
elements. Then $\mathcal{S}_n$ acts on $G^n$ by setting
\[
    s (\mathbf{g})
    =
    \left( g_{s^{-1}(1)}, \dots, g_{s^{-1}(n)}\right)
\]
for each $s \in \mathcal{S}_n$ and $\mathbf{g} = (g_1, \ldots, g_n)
\in G^n$. Then the wreath product $G(\mathcal{S}_n)$ of $G$ by
$\mathcal{S}_n$ is the semidirect product of $G^n$ by this action.
In particular, we have
\[
    ({\mathbf g},s)({\mathbf h}, t)
    =
    ({\mathbf g}\, s({\mathbf h}), st)
\]
and
\[
    ({\mathbf g},s)^{-1}
    =
    (s^{-1}({\mathbf g}^{-1}), s^{-1})
\]
for $\mathbf{g}, \mathbf{h} \in G^n$ and $s, t \in \mathcal{S}_n$.

\subsubsection{Conjugacy Classes}
\label{subsubsec-conjclasseswreath}

The conjugacy classes of $G(\mathcal{S}_n)$ correspond to conjugacy
classes of $G$ and conjugacy classes of $\mathcal{S}_n$.  First,
note that
\[
    ({\mathbf g},s) ({\mathbf h},t) ({\mathbf g},s)^{-1}
    =
    ({\mathbf g} s({\mathbf h}) sts^{-1} ({\mathbf g}^{-1}),\, sts^{-1})
\]
for $\mathbf{g}, \mathbf{h}\in G^n$ and $s,t \in \mathcal{S}_n$.
Hence, the conjugacy class of the element $({\mathbf h},s) \in
G(\mathcal{S}_n)$ is given by
\[
    \left\{ ({\mathbf g} s({\mathbf h}) sts^{-1}
    ({\mathbf g}^{-1}),\, sts^{-1}) | {\mathbf g} \in G^n ,  s\in \mathcal{S}_n\right\}.
\]
We can decompose each element of $G(\mathcal{S}_n)$ into a product
of cycles corresponding to cycles in $\mathcal{S}_n$; see \cite[page
129]{tamanoi1} or \cite[Section 1]{wang1}.  In fact, given
$({\mathbf g},s)\in G(\mathcal{S}_n)$, let $s = \prod_j s_j$ be the
disjoint cycle decomposition of $s \in \mathcal{S}_n$, unique up to
the order of the factors; note that we allow 1-cycles.  For each
$s_j=(j_1,\dots, j_r)$, let ${\mathbf g}_j = (g_{1,j}, \ldots ,
g_{n,j}) \in G^n$ denote the element whose $i$-component is equal to
$g_i$ if $i \in \{ j_1, \dots, j_r\}$ and $1$ otherwise.  For each
$j,k$, $({\mathbf g}_j, s_j)$ and $({\mathbf g}_k, s_k)$ commute,
and $({\mathbf g},s) = \prod_j ({\mathbf g}_j, s_j) \in
G(\mathcal{S}_n)$ is called the \emph{cycle decomposition} of
$({\mathbf g},s)$ corresponding to $s = \prod_j s_j$.

The following proposition, giving information about the conjugacy
class of $({\mathbf g},s)$ in $G(\mathcal{S}_n)$, is stated in
\cite[Proposition 3-1]{tamanoi1} for $G$ finite.   Notice that the
outside product in Equation \ref{eq-typedecomposition} is restricted
to a finite set; with this modification, the proof by direct
computation extends to the case of $G$ compact.  We let $G_\ast$
denote the set of conjugacy classes of $G$ and $(c)$ the conjugacy
class of $c \in G$.

\begin{proposition}
\label{prop-cycledecompwrprod}

Let $({\mathbf g},s)= \prod\limits_j ({\mathbf g}_j, s_j) \in
G(\mathcal{S}_n)$ be a disjoint cycle decomposition of a fixed
$({\mathbf g},s) \in G(\mathcal{S}_n)$ as above.

\begin{enumerate}

\item
Fix $j$, and let $r$ be the length of the cycle $s_j$.  We consider
$(\mathbf{g}_j, s_j)$ as an element of $\mathcal{S}_r \leq \mathcal{S}_n$ in
the usual way.  Let $c_j \in
G$ denote the \emph{cycle product} of $s_j$, i.e. the product of the
components of $\mathbf{g}_j$, so that $c_j = \prod_{i=1}^n g_{i,j}$.
Then if $\mathbf{d}_j = (g_{i_1}, g_{i_2} g_{i_1}, \dots, g_{i_{r}}
g_{i_{{r}-1}}\dots g_{i_1}) \in G^r$ and $\mathbf{c}_j =
(c_j,1,\dots,1) \in G^r$, we have
\[
    ({\mathbf g}_j, s_j) = ({\mathbf d }_j, 1) (\mathbf{c}_j, s_j) ({\mathbf d }_j,1 )^{-1} \in G(\mathcal{S}_{r}).
\]
Moreover, $({\mathbf d }_j, 1)$ and $({\mathbf d }_k, 1)$ commute
for any $j, k$.

\item
Let $G_\ast(\mathbf{g}, s)$ denote the collection of conjugacy
classes $(c) \in G_\ast$ such that $c_j \in (c)$ for some $j$.  For
each $(c) \in G_\ast(\mathbf{g}, s)$ and each $r = 1, 2, \ldots, n$,
let $m_r(c)$ denote the number of $r$-cycles $s_j$ of $s$ such that
the corresponding $c_j$ is conjugate to $c$.  We enumerate the
corresponding $s_j$ and $\mathbf{c}_j$ as $s_{j((c),r)_1}, \ldots
s_{j((c),r)_{m_r(c)}}$ and $\mathbf{c}_{j_1}, \ldots ,
\mathbf{c}_{j_{m_r(c)}}$, which are now considered as elements
of $\mathcal{S}_n$ and $G(\mathcal{S}_n)$, respectively. Let $({\mathbf d }, 1)=
\prod_j ({\mathbf d }_j,1) \in G(\mathcal{S}_r)$. Then the cycle
decomposition $({\mathbf g},s)= \prod_j ({\mathbf g}_j, s_j)$ of
$({\mathbf g},s)$ induces the decomposition
\begin{equation}
\label{eq-typedecomposition}
    ({\mathbf g},s)
    =
    ({\mathbf d },1)
    \left(
    \prod_{(c) \in G_\ast(\mathbf{g}, s)}
    \;\;
    \prod_{r = 1}^n
    \prod_{i=1}^{m_r(c)}
    ({\mathbf c}_{j((c),r)_i}, s_{j((c),r)_i})
    \right)
    ({\mathbf d },1)^{-1}
    \in G(\mathcal{S}_n).
\end{equation}
Note that $G_\ast(\mathbf{g}, s)$ is finite so that the above
product is defined.

\end{enumerate}

\end{proposition}

\begin{proof}

Part (1) can be verified by a direct computation.  Then (2) follows
as we are taking the product over all possible cycle lengths, and as
the $({\mathbf g}_j, s_j)$ commute.

\end{proof}

It follows from Proposition \ref{prop-cycledecompwrprod} that for
each cycle $s_j$, the conjugacy class of the corresponding $c_j$ is
uniquely determined.  Hence, we have the following.

\begin{definition}
\label{def-partitionvalfunction}

Fix $({\mathbf g},s) \in G(\mathcal{S}_n)$.  With the notation as in
Proposition \ref{prop-cycledecompwrprod}, there is a
partition-valued function $\rho : G_\ast \to \mathcal{P}$, where
$\mathcal{P}$ denotes the set of partitions of all nonnegative
integers.  By a partition of $m$, we mean a string $M_r$ indexed by
$r = 1, \ldots , m$ such that $\sum_{r=1}^m rM_r = m$. Indeed, every
$(c) \in G_\ast(\mathbf{g}, s)$ determines the partition $m_r(c)$,
$r=1, \dots, t_c$, of $\sum_{r} r m_r(c) = t_c$.  If $(c) \in
G_\ast, (c) \notin G_\ast(\mathbf{g}, s)$, then we simply set
$\rho(c) = 0$ the empty partition of $0$.

Note that for a given $(\mathbf{g}, s) \in G(\mathcal{S}_n)$ and
$(c) \in G_\ast$, there are only finitely many nonzero $m_r(c)$.
Moreover, as $s \in \mathcal{S}_n$ is decomposed into cycles,
\begin{equation}
\label{eq-partitionssumton}
    \sum\limits_{r=1}^n \sum\limits_{(c) \in G_\ast} r m_r(c) =   n.
\end{equation}
In particular, $t_c \leq n$ is finite.  If we set $m_r = \sum_{(c)
\in G_\ast} m_r(c)$, then $\sum_{r=1}^n rm_r = n$ so that the $m_r$
determine a partition of $n$.

The partition-valued function $\rho : G_\ast \rightarrow
\mathcal{P}$ associated to $({\mathbf g},s)$ is called the
\emph{type of $({\mathbf g},s)$}, and is formally denoted using the
notation
\[
    \rho(c)= (1^{m_1(c)}, \dots, r^{m_r(c)}, \ldots ).
\]
where $m_r(c) \not= 0$ only for finitely many $(c) \in G_\ast, r \in
\{1, \dots, n \}$.

\end{definition}

By Proposition \ref{prop-cycledecompwrprod}, the conjugacy classes
of $G(\mathcal{S}_n)$ are parameterized by types, generalizing to
wreath products the well-known fact that the conjugacy classes of
$\mathcal{S}_n$ are determined by the structure of their cycle
decomposition. For further generalizations to actions of wreath
products, see Section \ref{subsec-actionswreathprodmanifolds} and in
particular Remark \ref{rem-finiteconjclasses}.

\subsubsection{Centralizers}
\label{subsubsec-centralizerswreathrod}

We will now study the structure of the centralizer of an element
$({\mathbf g},s)$ of the wreath product $G(\mathcal{S}_n)$ up to
conjugacy by $G(\mathcal{S}_n)$. We will prove in Theorem
\ref{thrm-centeralizerwreathprod} that
$C_{G(\mathcal{S}_n)}({\mathbf g},s)$ can be decomposed as a product
of wreath products.

Let $({\mathbf g},s)= \prod_j ({\mathbf g}_j, s_j) \in
G(\mathcal{S}_n)$ be the cycle decomposition of $({\mathbf g},s)$.
If $({\mathbf h},t) \in G(\mathcal{S}_n)$ commutes with $({\mathbf
g},s)$, then
\[
    ({\mathbf h},t) ({\mathbf g},s) ({\mathbf h},t)^{-1} =({\mathbf g},s).
\]
It follows from Proposition \ref{prop-cycledecompwrprod} and the
fact that conjugation preserves cycles, cycle length, and type that
for each $j$, there is a unique $k$ such that
\[
    ({\mathbf h},t) ({\mathbf g}_j,s_j) ({\mathbf h},t)^{-1} =({\mathbf g}_k, s_k).
\]
Moreover, $(c_j) = (c_k) = (c) \in G_\ast(\mathbf{g}, s)$, and $s_i$
and $s_j$ are both of length $r$. It follows that conjugation by
$({\mathbf h},t)$ permutes the $m_r(c)$ elements $({\mathbf c
}_{j((c), r)_i}, s_{j((c), r)_i})$ in the decomposition given by
Equation \ref{eq-typedecomposition} leaving $r$ and $(c)$ fixed.  It
therefore induces a homomorphism
\[
    p: C_{G(\mathcal{S}_n)} ({\mathbf g},s)
    \rightarrow
    \prod_{(c)\in G_\ast(\mathbf{g}, s)}
    \;\;
    \prod_{r = 1}^n \mathcal{S}_{m_r(c)},
\]
where $\mathcal{S}_{m_r(c)}$ as usual denotes the permutation group
on $m_r(c)$ elements.  We have the following, stated as \cite[Lemmas
3-3 and 3-4]{tamanoi1} for $G$ finite.

\begin{proposition}
\label{prop-homsplitsurjectiveandcentralizer}

The homomorphism
\[
    p: C_{G(\mathcal{S}_n)} ({\mathbf g},s)
    \rightarrow
    \prod_{(c)\in G_\ast(\mathbf{g}, s)}
    \;\;
    \prod_{r = 1}^n \mathcal{S}_{m_r(c)}
\]
is split surjective, and
\begin{equation}
\label{eq-kerp}
    \mbox{Ker}(p)
    =
    \left\{ ({\mathbf h}, t) \in G(\mathcal{S}_n) : ({\mathbf h}, t)
    ({\mathbf g}_j, s_j) ({\mathbf h}, t)^{-1} = ({\mathbf g}_j, s_j) \; \forall j \right\}.
\end{equation}
If $G^r_{(j)}$ denotes the subgroup of $G(\mathcal{S}_n)$ isomorphic
to $G(\mathcal{S}_r)$ corresponding to the positions of $s_j$, then
\begin{equation}
\label{eq-wreathcentralizer}
    C_{G^r_{(j)}}({\mathbf g}_j, s_j)
    \cong
    C_G(c_j)\langle a_{r,c_j}\rangle
\end{equation}
where $a_{r,c_j} =({\mathbf c}_j, (12\dots r))= ( (c_j,1,\dots, 1),
(12\dots r))$.

\end{proposition}

By $C_G(c_j)\langle a_{r,c_j}\rangle$, we mean the subgroup of
$G(\mathcal{S}_r)$ generated by $C_G(c_j)$ and $a_{r,c_j}$. Note
that $(a_{r,c_j})^r = ((c_j, \ldots c_j), 1)$ and $[a_{r,c_j},
C_G(c_j)] = 1$, with $C_G(c_j)$ identified with a subgroup of
$G_{(j)}^r$ as diagonal elements $((h, \ldots, h), 1)$.

\begin{proof}

Following \cite[Lemma 3-3]{tamanoi1}, we write each $s_j$ in the
decomposition $({\mathbf g},s) =\prod_j ({\mathbf g}_j, s_j)$
starting with the smallest integer.  We construct a homomorphism
\[
    \lambda : \prod\limits_{(c)\in G_\ast(\mathbf{g}, s)} \;\;
    \prod\limits_{r=1}^n {\mathcal S}_{m_r(c)} \longrightarrow {\mathcal S}_n
\]
by sending an element $\overline{t} \in \prod_{(c)\in
G_\ast(\mathbf{g}, s)} \prod_{r=1}^n {\mathcal S}_{m_r(c)}$ that
sends the cycle $s_j = (j_1, \dots, j_r)$ to $s_k = (k_1, \dots,
k_r)$ to $\lambda (\overline{t})= t \in {\mathcal S}_n$, where
$t(j_\ell)= k_\ell$, $1 \leq \ell \leq r$. Then  each element of the
image of $\lambda$ commutes with $s$.  Define the homomorphism
\[
\begin{array}{rccl}
    \Lambda :&
    \prod\limits_{(c)\in G_\ast(\mathbf{g}, s)} \;\; \prod\limits_{r=1}^n {\mathcal S}_{m_r(c)}
    &\longrightarrow&
    C_{G({\mathcal S}_n)} ({\mathbf g}, s)
            \\\\
    :&\overline{t}
    &\longmapsto&
    ({\mathbf d }, 1)  (1, \lambda(\overline{t})) ({\mathbf d },
    1)^{-1}
\end{array}
\]
with $(\mathbf{d}, 1)$ as in Proposition
\ref{prop-cycledecompwrprod}.  Then as $\Lambda (\overline{t})$
commutes with $({\mathbf g}, s)$ for every $\overline{t} \in
\prod_{(c)\in G_\ast(\mathbf{g}, s)} \prod_{r=1}^n {\mathcal
S}_{m_r(c)}$ and $p\circ \Lambda$ is the identity, $p$ is split
surjective.  Equation \ref{eq-kerp} follows directly.

To prove Equation \ref{eq-wreathcentralizer}, first note that if
$({\mathbf h}, 1)$ commutes with $( {\mathbf c}_r, (12\dots r))$ for
some ${\mathbf h} \in G^r$, then ${\mathbf h}$ must be of the form
$(h,h,\dots, h)$ for some $h \in C_G(c)$.  Similarly, if $({\mathbf
h}, \sigma )$ commutes with $({\mathbf c}_r, (12\dots r))$, then
$\sigma \in C_{{\mathcal S}_r}(12\dots r)$. The result follows after
noting that by Proposition \ref{prop-cycledecompwrprod}, every
element of a given type can be conjugated to a product of elements
of the above type.

\end{proof}

With this, the following important description of the centralizer of
$({\mathbf g },s)$ in $G(\mathcal{S}_n)$ generalizes to the case of
$G$ compact, connected with identical proof; see \cite[Theorem
3-5]{tamanoi1} and \cite[Lemma 2]{wang1}.

\begin{theorem}
\label{thrm-centeralizerwreathprod}

Let $({\mathbf g }, s )$ in $G(\mathcal{S}_n)$.  Then
\[
    C_{G(\mathcal{S}_n) } ({\mathbf g}, s )
    \cong
    \prod_{(c) \in G_\ast(\mathbf{g}, s)}
    \;\;
    \prod_{r=1}^n
    C_G(c)  \langle a_{r,c}\rangle (\mathcal{S}_{m_r(c)}),
\]
where $C_G(c) \langle a_{r,c} \rangle (\mathcal{S}_{m_r(c)})$ is a
wreath product, $a_{r,c}^r = (c, \ldots c)$, and $m_r(c)$ again
denotes the number of $r$-cycles in the cycle decomposition of
$(\mathbf{g}, s)$ with cycle product in the conjugacy class $(c)$.
The isomorphism is induced by the conjugation given in Equation
\ref{eq-typedecomposition}.

\end{theorem}

Hence, the centralizer of an element $(\mathbf{g}, s) \in
G(\mathcal{S}_n)$ can be decomposed into a product of centralizers
according to the cycle decomposition of $s \in \mathcal{S}_n$.


\subsection{Actions of Wreath Products on Manifolds}
\label{subsec-actionswreathprodmanifolds}

Let $M$ be a $G$-manifold, and consider the action of $G({\mathcal
S}_n)$ on $M^n$ given by
\[
    ({\mathbf g},s)  (x_1, \dots, x_n)
    =
    (g_1 x_{s^{-1}(1)}, \dots, g_n x_{s^{-1}(n)}).
\]
We let $MG({\mathcal S}_n)$ denote the orbifold presented by the
groupoid $M^n \rtimes G({\mathcal S}_n)$, called the \emph{wreath
symmetric product of $M \rtimes G$}. When $n=1$, note that $M
\rtimes G(\mathcal{S}_1) = M \rtimes G$.

In this section, we study the structure of the fixed point sets
$(M^n)^{\langle (\mathbf{g}, s)\rangle}$ for a fixed $(\mathbf{g},
s) \in G({\mathcal S}_n)$.  The results of this section generalize
to almost free, effective actions of a compact, connected  Lie group
several results proven in \cite{wang1} and \cite{tamanoi1} for $G$
finite.  First, we note the following general facts about group
actions and symmetric products.

Suppose $K$ is a compact, connected Lie group acting almost freely,
effectively, and smoothly on the smooth manifold $Z$ in such a way
that $Z \rtimes K$ is Morita equivalent to an orbifold groupoid. By
\cite[Theorem 5.3]{ademruan} (see also \cite{farsiseaton2}),
\[
    H_\Z^\ast(Z\rtimes K)
    =
    \bigoplus_{(k) \in t_{Z;K}} H^\ast\left(Z^{\langle k\rangle}\rtimes C_K(k)\right),
\]
where $t_{Z;K}$ denotes the set of conjugacy classes in $K$ with
non-empty fixed point set in $Z$.  It follows that the delocalized
cohomology is isomorphic to the complexified orbifold $K$-theory in
this case.

If $K$ and $L$ are compact, connected Lie groups such that $L$ acts
on $K$ and $K \rtimes L$ acts on $Z$, then
\begin{equation}
\label{eq-doublecrossprod}
\begin{array}{rcl}
    H^\ast(Z\rtimes (K \rtimes L))
    &\cong&
    H^\ast(Z)^{K \rtimes L}
            \\\\
    &\cong&
    H^\ast(Z \rtimes L)^L
\end{array}
\end{equation}
where the latter two expressions indicate the invariant cohomology.

Now, consider the $\mathcal{S}_n$-action on $M^n$, $M$ a smooth
manifold, by permuting factors, and let $\sigma \in \mathcal{S}_n$.
Then we have that
\[
    C_{{\mathcal S}_n}(\sigma)
    \cong
    \mathcal{S}_{m_1(\sigma)} \times
    \prod\limits_{j=2}^n
    {\mathbb Z}_{m_j(\sigma)} \rtimes {\mathcal S}_{m_j(\sigma)}.
\]
Letting $\Delta_j(M)$ denote the diagonal in $M^j$, we have
\[
\begin{array}{rcl}
    (M^n)^{\langle\sigma\rangle}
    &\cong&
    \prod\limits_{j=1}^{n} (\Delta_j(M))^{m_{j}(\sigma)}
                \\\\
    &\cong&
    \prod\limits_{j=1}^n  M^{m_{j}(\sigma)}.
\end{array}
\]
Hence,
\[
    (M^n)^{\langle\sigma\rangle} \rtimes C_{{\mathcal S}_n}(\sigma )
    \cong
    \prod\limits_{j=1}^n  M^{m_j(\sigma)} \rtimes {\mathcal S}_{m_j(\sigma)}.
\]
By definition, the delocalized cohomology of $M^n \rtimes \mathcal{S}_n$ is given by
\begin{equation}
\label{eq-symmproduct}
    H^\ast \left( \widetilde{(M^n \rtimes \mathcal{S}_n)}_\Z \right)
    =
    \bigoplus_{(\sigma)\in t_{M^n; {\mathcal S}_n }}
    H^\ast(M^{\langle \sigma \rangle} \rtimes C_{{\mathcal S}_n}({\sigma})  ),
\end{equation}
where $t_{M^n; {\mathcal S}_n}$ again denotes the conjugacy classes
$\sigma \in \mathcal{S}_n$ such that $(M^n)^{\langle \sigma \rangle}
\neq \emptyset$.  As $\sigma$ is conjugate to $\tau$ in
$\mathcal{S}_n$ if and only if $m_j(\sigma) = m_j(\tau)$ for each $j
= 1,\ldots n$, the conjugacy classes of the symmetric group are
parametrized by the partitions of $n$. Summing over all $j, m_j \in
\{ 1, \ldots , n \}$ with $\sum_{j=1}^n jm_j = n$, Equation
\ref{eq-symmproduct} becomes
\[
    H^\ast \left( \widetilde{(M^n \rtimes \mathcal{S}_n)}_\Z \right)
    \cong
    \bigoplus_{j, m_j} H^\ast \left(
    \prod\limits_{j =1}^n M^{m_j}/ {\mathcal S}_{m_j} \right).
\]
Using Equation \ref{eq-doublecrossprod} and the K\"{u}nneth Theorem (see \cite{rosenbergschotet}),
it follows that
\[
\begin{array}{rcl}
    H^\ast \left( \widetilde{(M^n \rtimes \mathcal{S}_n)}_\Z \right)
    &\cong&
    \bigoplus\limits_{j,m_j} \;\; \bigotimes\limits_{j=1}^n
    \left(H^\ast(M)^{\bigoplus^{m_j}}\right)^{{\mathcal S}_{m_j}}
            \\\\
    &\cong&
    \bigotimes\limits_{j,m_j} \;\; \bigotimes\limits_{j=1}^n
    {\mathcal {SP}}^{m_j} [H^\ast(M)]
\end{array}
\]
where $\mathcal{SP}$ denote the symmetric product, i.e. the subspace fixed under the
action of the symmetric group by permuting coordinates.

We now return to the case of a wreath symmetric product.  We begin
by studying the special case of $(\mathbf{g}, s) = ((g, 1, \ldots,
1), \tau) \in G(\mathcal{S}_n)$ with $\tau = (1, 2, \ldots , n) \in
\mathcal{S}_n$ that arose in Proposition
\ref{prop-homsplitsurjectiveandcentralizer}.  Given Propositions
\ref{prop-cycledecompwrprod} and
\ref{prop-homsplitsurjectiveandcentralizer}, the proofs of
\cite[Lemmas 4 and 5]{wang1} generalize to the case of $G$ compact
and connected, yielding the following.

\begin{proposition}
\label{prop-wrprodfixedptspecific}

Suppose that $a = ((g, 1,\dots, 1), \tau) \in G({\mathcal S}_n)$ for
some $g \in G$ with $\tau = (1,2,\dots, n) \in {\mathcal S}_n$. Then
\[
    (M^n)^{\langle a \rangle}
    =
    \{ (x,\dots, x) \in M^n : x = gx\}.
\]
The group $C_{G({\mathcal S}_n)}(a)$ is isomorphic to $C_G(g)\langle
a\rangle$, and the action of  $C_G(g)\langle a\rangle$ on
$(M^n)^{\langle a\rangle}$ can be identified via the diagonal map
with the action of $C_G(g)$ on $M^{\langle g\rangle}$ together with
the trivial $a$-action. Therefore, by Equation
\ref{eq-doublecrossprod}, we have
\[
    H^\ast\left((M^n)^{\langle a\rangle}\rtimes C_{G({\mathcal S}_n)}(a)\right)
    \cong H^\ast\left(M^{\langle g\rangle}\right)^{C_G(g)}.
\]

\end{proposition}

Using the decomposition given in Proposition
\ref{prop-cycledecompwrprod}, it follows that $G(\mathcal{S}_n)$
acts locally freely on $M^n$ if and only if $G$ acts locally freely
on $M$.  Hence, $M \rtimes G$ is Morita equivalent to an orbifold
groupoid if and only if the same holds for $M^n \rtimes G({\mathcal
S}_n)$.  With this, we can determine the delocalized equivariant
cohomology of $MG({\mathcal S}_n)$, that is, the standard cohomology
of the inertia orbifold $\widetilde{MG({\mathcal S}_n)}_\Z$ of
$MG({\mathcal S}_n)$ in the case that $M \rtimes G$ is Morita
equivalent to an orbifold groupoid (see Subsection
\ref{subsec-cohomology}).

\begin{proposition}
\label{prop-wrprodfixedptgeneral}

For a fixed $({\mathbf g}, s ) \in G({\mathcal S}_n)$ with nonempty
fixed-point set in $M^n$, the groupoid $(M^n)^{\langle ({\mathbf g},
s ) \rangle}\rtimes C_{G({\mathcal S}_n)}(({\mathbf g}, s ))$ is
given as follows. First, using the notation of Section
\ref{subsec-wreathprodcompact}, let
\[
    Y^{m_r(c)} = \left(M^{\langle c \rangle}\right)^{m_r(c)}
    \rtimes
    (C_G(c))^{m_r(c)},
\]
and let ${\mathcal S}_{m_r(c)}$ act on $Y^{m_r(c)} $ by permuting
the coordinates of $(M^{\langle c\rangle})^{m_r(c)}$. Then
\[
    Y^{m_r(c)}\rtimes {\mathcal S}_{m_r(c)}
    = \left((M^{\langle c\rangle})^{m_r(c)} \rtimes C_G(c)^{m_r(c)} \right)
    \rtimes {\mathcal S}_{m_r(c)}.
\]
It follows by Equation \ref{eq-doublecrossprod} that
\[
    H^\ast \left(Y^{m_r(c)}\rtimes {\mathcal S}_{m_r(c)}\right)
    \cong H^\ast \left((M^{\langle c\rangle})^{m_r(c)} \rtimes C_G(c)^{m_r(c)}\right)^{{\mathcal
S}_{m_r(c)}}.
\]
Hence,
\[
    H^\ast\left((M^n)^{\langle a\rangle}\rtimes C_{G({\mathcal S}_n)}(a)\right)
    \cong
    \bigotimes\limits_{(c) \in t_{M; G}} \;\; \bigotimes\limits_{r=1}^n
    H^\ast
    \left((M^{\langle c \rangle})^{m_r(c)} \rtimes C_G(c)^{m_r(c)}\right)^{{\mathcal S}_{m_r(c)}},
\]
where $\sum_r r {m_r}(c)=t_c \leq n$.

\end{proposition}

\begin{proof}

The decomposition $({\mathbf g},s)= \prod_{i} ({\mathbf g}_i, s_i)$
in $G({\mathcal S}_n)$ from Proposition \ref{prop-cycledecompwrprod}
implies that
\[
    (M^n)^{\langle ( {\mathbf g},s) \rangle }
    \cong
    \prod\limits_j (M^{r_j})^{\langle( {\mathbf g}_j, s_j)\rangle},
\]
where $M^{r_j} \subseteq M^n$ corresponds to the positions of $s_j$.
By Proposition \ref{prop-cycledecompwrprod},
\[
    (M^n)^{\langle({\mathbf g},s)\rangle}
    \cong
    \prod\limits_{(c)\in G_\ast(\mathbf{g}, s)}
    \left(M^{\langle c \rangle}\right)^{{\sum}_{r=1}^n m_r (c)}.
\]
Note in particular that each cycle $s_j$ corresponds to a cycle product with conjugacy class
$(c)$ such that $M^{\langle c \rangle} \neq \emptyset$.

By synthesizing Theorem \ref{thrm-centeralizerwreathprod} and
Proposition \ref{prop-wrprodfixedptgeneral}, the action of the
centralizer $C_{G({\mathcal S}_n)}(\mathbf{g}, s)$ on
$(M^n)^{\langle (\mathbf{g}, s) \rangle}$ can be identified with a
product of actions of wreath products of the form
\[
    \left[ C_G(c) \langle a_{r, c}\rangle \right]^{m_r(c)}
    \rtimes {\mathcal S}_{m_r(c)},
\]
for an element $a_{r, c}$ acting trivially on
$M^{\langle c \rangle }$, and with ${\mathcal S}_{m_r(c)}$ acting by permutations
on $[C_{G(c)}a_{r, c}]^{m_r(c)}$.  The result follows.

\end{proof}

\begin{remark}
\label{rem-finiteconjclasses}

The above results imply that the conjugacy classes of
$G(\mathcal{S}_n)$ with non-empty fixed point sets in $M^n$ are
parameterized by the \emph{finite} set $t_{M;G}$ of conjugacy
classes $(c) \in G_\ast$ with non-empty fixed point sets and their
types $m_r(c)$ with $r = 1, \ldots, n$.  In particular, let $t_{M;G}
= \{ (c_1), \ldots , (c_N) \}$.  The set $t_{M^n; G(\mathcal{S}_n)}
= \{ (a) \in G(\mathcal{S}_n): (M^n)^{\langle a \rangle} \neq
\emptyset \}$ is parameterized by the $m_r(c_k)$ with $r = 1, \ldots
, n$ and $k = 1, \ldots N$. We emphasize the equalities $\sum_{(c)
\in t_{M; G}}\sum_r rm_r(c) = n$ and $\sum_r rm_r = n$ from
Definition \ref{def-partitionvalfunction}, which play an important
role in the sequel.

\end{remark}

By the above observations, when $M/G$ and hence $MG(\mathcal{S}_n)$
is an orbifold, the inertia orbifold
$\widetilde{(MG(\mathcal{S}_n))}_\Z$ of $MG(\mathcal{S}_n)$ is
presented by the groupoid
\begin{equation}
\label{eq-wreathinertiagroupoid}
\begin{array}{rcl}
    \mathbb{IMG}_n
    &=&
    \coprod\limits_{((\mathbf{g}, s)) \in t_{M^n; G(\mathcal{S}_n)}}
    (M^n)^{\langle (\mathbf{g}, s) \rangle} \rtimes C_{G(\mathcal{S}_n)}(\mathbf{g}, s)
                \\\\
    &=&
    \coprod\limits_{m_r(c)} \;
    \prod\limits_{r=1}^n \;
    \prod\limits_{(c) \in t_{M; G}}\;
    \left((M^{\langle c \rangle})^{m_r(c)}
    \rtimes (C_G(c)\langle a_{r,c} \rangle)^{m_r(c)}\right)
    \rtimes {\mathcal S}_{m_r(c)}.
\end{array}
\end{equation}
It follows that
\begin{equation}
\label{eq-wreathinertiacohom}
\begin{array}{rcl}
    H^\ast\left(\widetilde{(MG(\mathcal{S}_n))}_\Z\right)
    &=&
    \bigoplus\limits_{((\mathbf{g}, s)) \in t_{M^n; G(\mathcal{S}_n)}}
    H^\ast\left((M^n)^{\langle (\mathbf{g}, s) \rangle}
    \rtimes C_{G({\mathcal S}_n)}(\mathbf{g}, s)\right)
            \\\\
    &\cong&
    \bigoplus\limits_{m_r(c)} \;
    \bigotimes\limits_{r=1}^n \;
    \bigotimes\limits_{(c) \in t_{M; G}}
    H^\ast \left((M^{\langle c\rangle})^{m_r(c)}\rtimes (C_G(c)^{m_r(c)}
    \rtimes {\mathcal S}_{m_r(c)})\right)
            \\\\
    &\cong&
    \bigoplus\limits_{m_r(c)}\;
    \bigotimes\limits_{r=1}^n\;
    \bigotimes\limits_{(c) \in t_{M; G}}
    H^\ast \left((M^{\langle c\rangle})^{m_r(c)} \rtimes C_G(c)^{m_r(c)}\right)^{{\mathcal
    S}_{m_r(c)}}.
\end{array}
\end{equation}


\section{Behavior of $\Gamma$-Sectors Under Operations on Groups and Orbifolds}
\label{sec-operationsgroupsorbifolds}

In this section, examine operations on groups and orbifolds and the
corresponding operations on the $\Gamma$-sectors.


\subsection{The $\Gamma_1\times\Gamma_2$-Sectors}
\label{subsec-GammatimesGamma}

Let $\Gamma_1$ and $\Gamma_2$ be finitely generated groups.  Since
$\tilde{Q}_{\Gamma_1}$ is an orbifold with orbifold structure given
by the groupoid $\mathcal{G}^\Gamma = \mathcal{G} \ltimes
\mathcal{S}_{\mathcal G}^{\Gamma_1}$, it makes sense to form
the $\Gamma_2$-sectors of this orbifold as the
translation groupoid  $(\mathcal{G} \ltimes
\mathcal{S}_\mathcal{G}^{\Gamma_1}) \ltimes \mathcal{S}_{\mathcal{G}
\ltimes \mathcal{S}_{\mathcal G}^{\Gamma_1}}^{\Gamma_2}$. We denote
the orbifold associated to the resulting action groupoid
$\widetilde{\left(\tilde{Q}_{\Gamma_1}\right)}_{\Gamma_2}$. We claim
the following, which by virtue of \cite[Theorems 3.1 and 3.5]{farsiseaton2}
is a generalization of \cite[Proposition 2-1 (2)]{tamanoi2} (see
also \cite[Proposition 2-1]{tamanoi1}).

\begin{theorem}
\label{thrm-Gamma1timesGamma2formula}

Let $\mathcal{G}$ be a groupoid and $\Gamma_1$ and $\Gamma_2$
finitely generated discrete groups.  There is a groupoid isomorphism from
$\left(\mathcal{G} \ltimes \mathcal{ S}_{\mathcal
G}^{\Gamma_1}\right) \ltimes {\mathcal S}_{\mathcal{G} \ltimes
\mathcal{S}_\mathcal{G}^{\Gamma_1}}^{\Gamma_2}$ to $\mathcal{G}
\ltimes \mathcal{ S}_{\mathcal G}^{\Gamma_1 \times \Gamma_2}$. In
particular, if $\mathcal{G}$ is an orbifold groupoid, then the
orbifolds $\widetilde{\left(\tilde{Q}_{\Gamma_1}\right)}_{\Gamma_2}$
and $\tilde{Q}_{\Gamma_1 \times \Gamma_2}$, presented by
$\left(\mathcal{G} \ltimes \mathcal{ S}_{\mathcal
G}^{\Gamma_1}\right) \ltimes {\mathcal S}_{\mathcal{G} \ltimes
\mathcal{S}_\mathcal{G}^{\Gamma_1}}^{\Gamma_2}$ and
$\mathcal{G} \ltimes \mathcal{ S}_{\mathcal G}^{\Gamma_1 \times \Gamma_2}$,
respectively, are diffeomorphic.

\end{theorem}

\begin{proof}

First, we demonstrate a bijection between the spaces of objects. An
element of $\mathcal{ S}_{\mathcal G}^{\Gamma_1 \times \Gamma_2}$ is
a homomorphism $\Phi_x : \Gamma_1 \times \Gamma_2 \rightarrow G_x$
where $x \in G_0$.  An element of ${\mathcal S}_{\mathcal{G} \ltimes
\mathcal{S}_\mathcal{G}^{\Gamma_1}}^{\Gamma_2}$ is a homomorphism
$\psi_{\phi_x}$ from $\Gamma_2$ into the isotropy group of
the point $\phi_x \in \left(\mathcal{G} \ltimes \mathcal{S}_\mathcal{G}^{\Gamma_1}\right)_0 =
\mathcal{S}_\mathcal{G}^{\Gamma_1}$ with respect to the groupoid
$\mathcal{G} \ltimes \mathcal{S}_\mathcal{G}^{\Gamma_1}$. Note that the isotropy group of
$\phi_x$ is $C_{G_x}(\phi_x)$.

Fix $x \in G_0$; we claim that there is a bijection between
homomorphisms $\Phi_x : \Gamma_1 \times \Gamma_2 \rightarrow G_x$
and pairs of homomorphisms $\phi_x : \Gamma_1 \rightarrow G_x$ and
$\psi_{\phi_x} : \Gamma_2 \rightarrow C_{G_x}(\phi_x)$.  Let $\phi_x
: \Gamma_1 \rightarrow G_x$ and  $\psi_{\phi_x} : \Gamma_2
\rightarrow C_{G_x}(\phi_x)$. We define $\phi_x \centerdot
\psi_{\phi_x}$ to be the pointwise product, i.e.
\begin{equation}
\label{eq-inducdeftimeshoms}
\begin{array}{rccl}
    \phi_x \centerdot \psi_{\phi_x}  :& \Gamma_1 \times \Gamma_2    &  \longrightarrow & G_x
                            \\\\
                    :& (\gamma_1, \gamma_2)        &   \longmapsto     &
                    \phi_x(\gamma_1)\psi_{\phi_x}(\gamma_2) .
\end{array}
\end{equation}
As $\psi_{\phi_x}(\gamma_2) \in C_{G_x}(\phi)$ for each $\gamma_2 \in \Gamma_2$
so that it commutes with $\phi_x(\gamma_1)$ for each $\gamma_1 \in \Gamma_1$,
$\phi_x \centerdot \psi_{\phi_x}$ is clearly a homomorphism.

Given a homomorphism $\Phi_x : \Gamma_1 \times \Gamma_2 \rightarrow G_x$,
if we set $\phi_x(\gamma_1)  =   \Phi_x(\gamma_1, 1)$ for each $\gamma_1 \in \Gamma_1$ and
$\psi_{\phi_x}(\gamma_2)  =   \Phi_x(1, \gamma_2)$
for each $\gamma_2 \in \Gamma_2$, then $\Phi_x = \phi_x \centerdot \psi_{\phi_x}$.
Moreover, as
\[
\begin{array}{rcl}
    \phi_x(\gamma_1)\psi_{\phi_x}(\gamma_2)
        &=&     \Phi_x[(\gamma_1, 1)(1, \gamma_2)]                \\\\
        &=&     \Phi_x[(1, \gamma_2)(\gamma_1, 1)]                \\\\
        &=&     \psi_{\phi_x}(\gamma_2)\phi_x(\gamma_1),
\end{array}
\]
it is clear that $\psi_{\phi_x}(\gamma_2) \in C_{G_x}(\phi_x)$ for
each $\gamma_2 \in \Gamma_2$.

It follows that for each $x \in G_0$, the fiber of
\[
    \beta_{\Gamma_1 \times \Gamma_2}:
    \mathcal{S}_\mathcal{G}^{\Gamma_1 \times \Gamma_2}
    \longrightarrow
    G_0
\]
over each $x \in G_0$ is in bijective correspondence with the fiber of the composition
\[
    {\mathcal S}_{\mathcal{G} \ltimes \mathcal{S}_\mathcal{G}^{\Gamma_1}}^{\Gamma_2}
    \stackrel{\beta_{\Gamma_1}}{\longrightarrow}
    \mathcal{S}_\mathcal{G}^{\Gamma_1}
    \stackrel{\beta_{\Gamma_2}}{\longrightarrow}
    G_0
\]
over $x$.  Hence, the map
\[
\begin{array}{rccl}
    e_0  :& {\mathcal S}_{\mathcal{G} \ltimes \mathcal{S}_\mathcal{G}^{\Gamma_1}}^{\Gamma_2}
                &  \longrightarrow &        \mathcal{S}_\mathcal{G}^{\Gamma_1 \times \Gamma_2}
                                    \\\\
                    :& \psi_{\phi_x}        &   \longmapsto     &      \phi_x \centerdot
                    \psi_{\phi_x}
\end{array}
\]
is bijective.

Now, we represent $\left(\mathcal{G} \ltimes \mathcal{ S}_{\mathcal
G}^{\Gamma_1}\right) \ltimes {\mathcal S}_{\mathcal{G} \ltimes
\mathcal{S}_\mathcal{G}^{\Gamma_1}}^{\Gamma_2}$ as a
$\mathcal{G}$-space.  An arrow in $\mathcal{G} \ltimes
\mathcal{S}_\mathcal{G}^{\Gamma_1}$ is given by a pair $(h, \phi_x)$
with $\phi_x \in \mathcal{S}_\mathcal{G}^{\Gamma_1}$ and an $h \in
\mathcal{G}$ such that $s(h) = x$. Then for each $\psi_{\phi_x} \in
\mathcal{S}_{\mathcal{G}\ltimes
\mathcal{S}_\mathcal{G}^{\Gamma_1}}^{\Gamma_2}$ and $\gamma_2 \in
\Gamma_2$, the action is given by
\[
    [(h, \phi_x)\cdot \psi_{\phi_x}] (\gamma_2)
    =
    h \psi_{\phi_x}(\gamma_2) h^{-1}
\]
where $h \psi_{\phi_x}(\gamma_2) h^{-1}$ is a homomorphism from
$\Gamma_2$ to $C_{G_{t(h)}}(h\phi_x h^{-1})$. Hence, if we let $\alpha
: \mathcal{S}_{\mathcal{G}\ltimes
\mathcal{S}_\mathcal{G}^{\Gamma_1}}^{\Gamma_2} \rightarrow G_0$ be
defined by $\alpha(\psi_{\phi_x}) = x$, then $\alpha$ is the anchor
map of a $\mathcal{G}$-action defined by
\[
    [h \cdot \psi_{\phi_x}] (\gamma_2)
    =
    h \psi_{\phi_x}(\gamma_2) h^{-1},
\]
defined whenever $\alpha(\psi_{\phi_x}) = s(h) = x$.  The
requirements of a groupoid action are clearly satisfied, so that
$\left(\mathcal{G} \ltimes \mathcal{ S}_{\mathcal
G}^{\Gamma_1}\right) \ltimes {\mathcal S}_{\mathcal{G} \ltimes
\mathcal{S}_\mathcal{G}^{\Gamma_1}}^{\Gamma_2}$ is isomorphic to
$\mathcal{G} \ltimes {\mathcal S}_{\mathcal{G} \ltimes
\mathcal{S}_\mathcal{G}^{\Gamma_1}}^{\Gamma_2}$.

With this, note that for each $h \in G_1$ with $s(h) = x$,
$\psi_{\phi_x} \in \mathcal{S}_{\mathcal{G} \ltimes
\mathcal{S}_\mathcal{G}^{\Gamma_1}}^{\Gamma_2}$, $\gamma_1 \in
\Gamma_1$ and $\gamma_2 \in \Gamma_2$,
\[
\begin{array}{rcl}
    e_0[h \cdot \psi_{\phi_x}](\gamma_1, \gamma_2)
        &=&     e_0[h\psi_{\phi_x}h^{-1}](\gamma_1, \gamma_2)
                                                                \\\\
        &=&     [(h\cdot \phi_x) \centerdot (h \psi_{\phi_x}h^{-1})](\gamma_1, \gamma_2)             \\\\
        &=&     (h \phi_x(\gamma_1) h^{-1})(h \psi_{\phi_x}(\gamma_2) h^{-1})                                  \\\\
        &=&     h (\phi_x(\gamma_1)\psi_{\phi_x}(\gamma_2)) h^{-1}                                  \\\\
        &=&     [h\cdot (\phi_x \centerdot \psi_{\phi_x})](\gamma_1, \gamma_2)            \\\\
        &=&     [h \cdot e_0(\psi_{\phi_x})](\gamma_1, \gamma_2)
\end{array}
\]
It follows that $e_0$ is $\mathcal{G}$-equivariant, and hence by
\cite[Lemma 2.6]{farsiseaton2}, is
the map on objects of a groupoid isomorphism. The map on arrows is
given by
\[
\begin{array}{rccl}
    e_1  :& \mathcal{G} \ltimes {\mathcal S}_{\mathcal{G} \ltimes \mathcal{S}_\mathcal{G}^{\Gamma_1}}^{\Gamma_2}
                &  \longrightarrow &        \mathcal{G} \ltimes \mathcal{S}_\mathcal{G}^{\Gamma_1 \times \Gamma_2}
                                    \\\\
                    :& (h, \psi_{\phi_x})        &   \longmapsto     &      (h, \phi_x
                    \centerdot
                    \psi_{\phi_x}).
\end{array}
\]

Now, assume that $\mathcal{G}$ is an orbifold groupoid.  We claim
that $e_0$ is smooth, completing the proof.

Recall from \cite[Lemma 2.2]{farsiseaton1} the smooth structures on
${\mathcal S}_{\mathcal{G} \ltimes
\mathcal{S}_\mathcal{G}^{\Gamma_1}}^{\Gamma_2}$ and
$\mathcal{S}_\mathcal{G}^{\Gamma_1 \times \Gamma_2}$.  Near $\phi_x
\in \mathcal{S}_\mathcal{G}^{\Gamma_1}$, a manifold chart is of the
form $\{ V_x^{\langle \phi_x \rangle}, \kappa_{\phi_x} \}$ with $\kappa_{\phi_x}$
defined as follows.  For each
$y \in V_x^{\langle \phi_x \rangle}$, there is an injective group homomorphism
$\xi_x^y : G_y \rightarrow G_x$ arising from the isomorphism between
$\mathcal{G}|_{V_x}$ and $G_x \ltimes V_x$.  As $y$ is fixed by $\Imt \phi_x$,
the image of $\xi_x^y$ contains $\Imt \phi_x$, so we can define
$\kappa_{\phi_x}(y) = (\xi_x^y)^{-1} \circ \phi_x : \Gamma \rightarrow G_y$.

Near
$\psi_{\phi_x} \in {\mathcal S}_{\mathcal{G} \ltimes
\mathcal{S}_\mathcal{G}^{\Gamma_1}}^{\Gamma_2}$, then, a manifold
chart is of the form $\left\{ \left( V_x^{\langle \phi_x \rangle}
\right)^{\langle \psi_{\phi_x} \rangle}, \kappa_{\psi_{\phi_x}}
\right\}$ with $\left(V_x^{\langle \phi_x \rangle}
\right)^{\langle \psi_{\phi_x} \rangle} = V_x^{\langle \phi_x,
\psi_{\phi_x} \rangle}$.  It follows from the definition of
$\kappa_{\psi_{\phi_x}}$ and the fact that the arrows in
$\mathcal{G}\ltimes\mathcal{S}_\mathcal{G}^{\Gamma_1}$ are those of
$\mathcal{G}$ that
\begin{equation}
\label{eq-induckappa}
    \kappa_{\psi_{\phi_x}} : y \longmapsto
    (\xi_x^y)^{-1}\circ\psi_{\phi_x}
    \in\mbox{HOM}(\Gamma_2, C_{G_x}((\xi_x^y)^{-1}\circ\phi_x)).
\end{equation}

Near $\Phi_x \in \mathcal{S}_\mathcal{G}^{\Gamma_1 \times
\Gamma_2}$, a manifold chart is of the form $\{ V_x^{\langle \Phi_x
\rangle}, \kappa_{\Phi_x} \}$. We have
\[
    \kappa_{\Phi_x} : y \longmapsto (\xi_x^y)^{-1} \circ \Phi_x \in
    \mbox{HOM}(\Gamma_1 \times \Gamma_2, G_y).
\]
Note that
\[
\begin{array}{rcl}
    V_x^{\langle e_0(\psi_{\phi_x}) \rangle}
        &=&     V_x^{\langle \phi_x \centerdot \psi_{\phi_x} \rangle} \\\\
        &=&     V_x^{\langle \phi_x, \psi_{\phi_x} \rangle},
\end{array}
\]
so that a manifold chart at $e_0(\psi_{\phi_x})$ for ${\mathcal
S}_\mathcal{G}^{\Gamma_1 \times \Gamma_2}$ is identified with a
chart at $\psi_{\phi_x}$ for ${\mathcal S}_{\mathcal{G} \ltimes
\mathcal{S}_\mathcal{G}^{\Gamma_1}}^{\Gamma_2}$.

Applying $e_0$ to Equation \ref{eq-induckappa}, we have that
\[
\begin{array}{rcl}
    e_0[\kappa_{\psi_{\phi_x}}(y)]
        &=&     e_0[(\xi_x^y)^{-1}\circ\psi_{\phi_x}]
                                                            \\\\
        &=&     [((\xi_x^y)^{-1}\circ\phi_x)
                \centerdot((\xi_x^y)^{-1}\circ\psi_{\phi_x})]
                                                            \\\\
        &=&     [(\xi_x^y)^{-1}\circ(\phi_x \centerdot \psi_{\phi_x})]
                                                            \\\\
        &=&     \kappa_{e_0(\psi_{\phi_x})}(y).
\end{array}
\]
It follows that on each chart, $e_0$ is simply the identity map on charts and
hence a diffeomorphism.  It follows from \cite[Lemma 2.6]{farsiseaton2}
that $e$ is an isomorphism of orbifold groupoids.

\end{proof}


\subsection{The $\Gamma$-Sectors of Product Orbifolds}
\label{subsec-productorbifolds}

In this subsection, we prove the following.

\begin{proposition}
\label{prop-productformula}

Suppose $\mathcal{G}$ and $\mathcal{H}$ are groupoids and $\Gamma$
is a finitely generated discrete group. Then $(\mathcal{G} \times
\mathcal{H})^\Gamma$ and $\mathcal{G}^\Gamma \times
\mathcal{H}^\Gamma$ are isomorphic. If $\mathcal{G}$ and
$\mathcal{H}$ are orbifold groupoids, then the isomorphism is of
orbifold groupoids.

\end{proposition}

Note that this in particular implies that if $Q_1$ is presented by $\mathcal{G}$
and $Q_2$ is presented by $\mathcal{H}$, then $\widetilde{(Q_1 \times
Q_2)}_\Gamma$ and $\widetilde{(Q_1)}_\Gamma \times
\widetilde{(Q_2)}_\Gamma$ are diffeomorphic orbifolds,
generalizing \cite[Proposition 2-1 (1)]{tamanoi2} to the case of
general orbifold groupoids (see also \cite[Proposition
2-1]{tamanoi1}).

\begin{proof}

Recall that $(\mathcal{G} \times \mathcal{H})^\Gamma =
(\mathcal{G}\times\mathcal{H})\ltimes(\mathcal{S}_\mathcal{G}^\Gamma\times\mathcal{S}_\mathcal{H}^\Gamma)$
and
$\mathcal{G}^\Gamma\times\mathcal{H}^\Gamma=(\mathcal{G}\ltimes\mathcal{S}_\mathcal{G}^\Gamma)
\times(\mathcal{H}\ltimes\mathcal{S}_\mathcal{H}^\Gamma)$. Note that
an element of $\mathcal{S}_\mathcal{G}^\Gamma \times
\mathcal{S}_\mathcal{H}^\Gamma$ is of the form $(\phi_x, \psi_y)$
where $x \in G_0$, $y \in H_0$, $\phi_x \in \mbox{HOM}(\Gamma,
G_x)$, and $\psi_y \in \mbox{HOM}(\Gamma, H_y)$.  Define the map
\[
\begin{array}{rcccl}
    e_0 :   &   \mathcal{S}_\mathcal{G}^\Gamma \times
            \mathcal{S}_\mathcal{H}^\Gamma
            &\longrightarrow
            &\mathcal{S}_{\mathcal{G} \times \mathcal{H}}^\Gamma              \\\\
            :&(\phi_x, \psi_y)  &\longmapsto&   \phi_x \times \psi_y
\end{array}
\]
where $\phi_x \times \psi_y \in \mbox{HOM}(\Gamma, G_x\times H_y)$
is defined by
\[
    (\phi_x \times \psi_y)(\gamma) = (\phi_x(\gamma), \psi_y (\gamma)) \in G_x \times H_y.
\]
That $e_0$ is bijective is obvious,
as it is inverted by $(pr_1, pr_2)$ where the $pr_i$ correspond to
projection onto the first and second factors, respectively.

Now, note that the componentwise $\mathcal{G}$- and
$\mathcal{H}$-actions defining
$\mathcal{G}^\Gamma\times\mathcal{H}^\Gamma$ correspond to commuting
$\mathcal{G}$- and $\mathcal{H}$-actions on
$\mathcal{S}_\mathcal{G}^\Gamma \times
\mathcal{S}_\mathcal{H}^\Gamma$.  In other words, if we define
$\mathcal{G}$- and $\mathcal{H}$-actions on
$\mathcal{S}_\mathcal{G}^\Gamma \times
\mathcal{S}_\mathcal{H}^\Gamma$ by $g\cdot(\phi_x, \psi_y) = (g\cdot
\phi_x, \psi_y)$ and $h\cdot(\phi_x, \psi_y) = (\phi_x, h\cdot
\psi_y)$ for $g \in G_1$, $h \in H_1$, and $(\phi_x, \psi_y) \in
\mathcal{S}_\mathcal{G}^\Gamma \times
\mathcal{S}_\mathcal{H}^\Gamma$, then we clearly have
\[
    g\cdot[h\cdot(\phi_x, \psi_y)]
        =       h\cdot[g\cdot(\phi_x, \psi_y)] .
\]
Moreover, with respect to this action,
\[
\begin{array}{rcl}
    \mathcal{G}^\Gamma\times\mathcal{H}^\Gamma
        &=& (\mathcal{G}\ltimes\mathcal{S}_\mathcal{G}^\Gamma) \times(\mathcal{H}\ltimes\mathcal{S}_\mathcal{H}^\Gamma)
                            \\\\
        &=& (\mathcal{G} \times \mathcal{H}) \ltimes
        (\mathcal{S}_\mathcal{G}^\Gamma \times
        \mathcal{S}_\mathcal{G}^\Gamma).
\end{array}
\]
Note further that for each $g\in G_1$, $h \in H_1$, and $(\phi_x,
\psi_y) \in \mathcal{S}_\mathcal{G}^\Gamma \times
\mathcal{S}_\mathcal{H}^\Gamma$, we have
\[
\begin{array}{rcl}
    e_0[ g \cdot[h\cdot (\phi_x, \psi_y)]]
        &=&     e_0(g \cdot \phi_x, h\cdot \psi_y)              \\\\
        &=&     (g\cdot \phi_x) \times (h\cdot \psi_y)          \\\\
        &=&     (g, h) \cdot (\phi_x \times \psi_y)                   \\\\
        &=&     (g, h) \cdot e_0 (\phi_x, \psi_y).
\end{array}
\]
Hence, the $\mathcal{G} \times \mathcal{H}$-actions on
$\mathcal{S}_\mathcal{G}^\Gamma \times
\mathcal{S}_\mathcal{H}^\Gamma$ and $\mathcal{S}_{\mathcal{G} \times
\mathcal{H}}^\Gamma$ obviously coincide via $e_0$, and by
\cite[Lemma 2.6]{farsiseaton2}, $e_0$ induces a
groupoid isomorphism.

If $\mathcal{G}$ and $\mathcal{H}$ are orbifold groupoids, then
restricted to charts of the form
\[
    \left\{ V_x^{\langle \phi_x\rangle}, C_{G_x}(\phi_x), \pi_x^{\phi_x} \right\} \times \left\{
    W_y^{\langle \psi_y \rangle}, C_{H_y}(\psi_y), \pi_y^{\psi_y} \right\}
\]
and
\[
    \left\{ (V \times W)_{(x,y)}^{\langle \phi_x \times \psi_y \rangle}, C_{G_x \times H_y}(\phi_x \times
    \psi_y), \pi_{(x, y)}^{\phi_x \times \psi_y} \right\},
\]
$e_0$ corresponds to the map $V_x^{\langle \phi_x \rangle} \times
W_y^{\langle \psi_y \rangle} \rightarrow (V_x \times W_y)^{\langle
\phi_x \times \psi_y \rangle}$ and hence is a diffeomorphism.  It
follows that $e$ is an isomorphism of orbifold groupoids.

\end{proof}

In \cite{farsiseaton1}, an equivalence relation
$\approx$ was defined on the $\mathcal{S}_\mathcal{G}^\Gamma$, and it was shown that
$\phi_x \approx \psi_y$ if and only if the $\mathcal{G}$-orbits of
$\phi_x$ and $\psi_y$ are on the same connected component of $|\mathcal{G}^\Gamma|$.
The $\approx$-class of $\phi_x$ is denoted $(\phi)_\approx$ or simply $(\phi)$
if no confusion is introduced.  Then $T_Q^\Gamma$ was defined to be the set of
$\approx$-classes of $\mathcal{S}_\mathcal{G}^\Gamma$.

As a consequence of Proposition \ref{prop-productformula},
the connected components of $\widetilde{(Q_1 \times Q_2)}_\Gamma$
clearly correspond to products of connected components of
$\widetilde{(Q_1)}_\Gamma$ and connected components of
$\widetilde{(Q_2)}_\Gamma$.  Hence, there is a bijection between
$T_{Q_1}^\Gamma \times T_{Q_2}^\Gamma$ and $T_{Q_1 \times
Q_2}^\Gamma$ as sets. This bijection is evidently given by
$((\phi)_{\approx, \mathcal{G}^\Gamma}, (\psi)_{\approx,
\mathcal{H}^\Gamma}) \mapsto (\phi \times \psi)_{\approx,
(\mathcal{G}\times\mathcal{H})^\Gamma}$.


\subsection{Maps on $\Gamma$-Sectors Induced by Group Homomorphisms}
\label{subsec-grouphoms}

Let $\Lambda$ and $\Gamma$ be finitely generated discrete groups and $\Phi :
\Lambda \rightarrow \Gamma$ a group homomorphism.  Given a groupoid
$\mathcal{G}$, $\Phi$ induces a map
\[
\begin{array}{rccl}
    e_{\Phi 0}     :&     \mathcal{S}_\mathcal{G}^\Gamma
                            &\longrightarrow& \mathcal{S}_\mathcal{G}^\Lambda
                                        \\\\
                   :&     \phi_x      &\longmapsto&   \phi_x \circ \Phi .
\end{array}
\]
As the notation suggests, this map is the map on objects of a
groupoid homomorphism.

\begin{lemma}
\label{lem-ephihom}

Let $\mathcal{G}$ be a groupoid, let $\Lambda$ and $\Gamma$ be finitely
generated discrete groups, and let $\Phi : \Lambda \rightarrow \Gamma$ be a group
homomorphism. Then $\Phi$ induces a groupoid homomorphism
\[
    e_\Phi : \mathcal{G}^\Gamma \longrightarrow \mathcal{G}^\Lambda .
\]
If $\mathcal{G}$ is an orbifold groupoid, then $e_\Phi$ is a
homomorphism of orbifold groupoids.  The map $e_{\Phi 0}$ on objects
is an immersion.

\end{lemma}

\begin{proof}

By \cite[Lemma 2.6]{farsiseaton2}, it suffices
to check that $e_0$ is $\mathcal{G}$-equivariant and, in the case of an orbifold groupoid,
smooth. Recall
that the anchor map of the $\mathcal{G}$-action on
$\mathcal{S}_\mathcal{G}^\Gamma$ is $\beta_\Gamma : \phi_x \mapsto
x$; similarly, the anchor map of the $\mathcal{G}$-action on
$\mathcal{S}_\mathcal{G}^\Lambda$ is $\beta_\Lambda : \psi_x \mapsto
x$.  Given $\phi_x \in \mathcal{S}_\mathcal{G}^\Gamma$, as $e_{\Phi
0}(\phi_x) = \phi_x \circ \Phi$ is a homomorphism from $\Lambda$
into $G_x$, we see that $\beta_\Lambda \circ e_0 (\phi_x) = x =
\beta_\Gamma(\phi_x)$, and so $\beta_\Lambda \circ e_{\Phi 0} =
\beta_\Gamma$.

For each $h \in G_1$, $\phi_x \in \mathcal{S}_\mathcal{G}^\Gamma$,
and $\gamma \in \Gamma$, we have
\[
\begin{array}{rcl}
    [h\cdot e_{\Phi 0}(\phi_x)](\gamma)
        &=&     h[\phi_x (\Phi(\gamma))]h^{-1}              \\\\
        &=&     (h\cdot \phi_x)(\Phi(\gamma))                     \\\\
        &=&     [e_{\Phi 0}(h\cdot\phi_x)](\gamma),
\end{array}
\]
so that $e_{\Phi 0}$ is a $\mathcal{G}$-equivariant map.

Now suppose $\mathcal{G}$ is an orbifold groupoid.  Fix $\phi_x \in
\mathcal{S}_\mathcal{G}^\Gamma$ and a linear chart $\{ V_x, G_x,
\pi_x \}$ for $Q$.  Then $\left\{ V_x^{\langle \phi_x \rangle},
\kappa_{\phi_x} \right\}$ is a manifold chart for
$\mathcal{S}_\mathcal{G}^\Gamma$ at $\phi_x$ and $\left\{
V_x^{\langle e_{\Phi 0}(\phi_x) \rangle}, \kappa_{e_{\Phi
0}(\phi_x)} \right\}$ is a manifold chart for
$\mathcal{S}_\mathcal{G}^\Lambda$.  Noting that $\Imt [e_{\Phi
0}(\phi_x)] = \Imt (\phi_x \circ \Phi) \leq \Imt \phi_x$ so that
$V_x^{\langle \phi_x \rangle} \subseteq V_x^{\langle e_{\Phi
0}(\phi_x) \rangle}$, we have that
\begin{equation}
\label{eq-ePhi0local}
    V_x^{\langle \phi_x \rangle}
    \stackrel{\kappa_{\phi_x}}{\longrightarrow}
    \mathcal{S}_\mathcal{G}^\Gamma
    \stackrel{e_{\Phi 0}}{\longrightarrow}
    \mathcal{S}_\mathcal{G}^\Lambda
    \stackrel{\kappa_{e_{\Phi 0}(\phi_x)}^{-1}}{\longrightarrow}
    V_x^{\langle e_{\Phi 0}(\phi_x) \rangle}
\end{equation}
is simply the embedding of $V_x^{\langle \phi_x \rangle}$ into
$V_x^{\langle e_{\Phi 0}(\phi_x) \rangle}$ as a subset and hence
smooth; in particular, $e_{\Phi 0}$ is an immersion.  Therefore,
$e_{\Phi 0}$ is a smooth, $\mathcal{G}$-equivariant map of
$\mathcal{S}_\mathcal{G}^\Gamma$ into
$\mathcal{S}_\mathcal{G}^\Lambda$, so that $e_\Phi$ is a Lie
groupoid homomorphism.

\end{proof}

We restrict our attention to the case of an orbifold groupoid
$\mathcal{G}$.  Using the identification of the multi-sectors with
the $\Gamma$-sectors corresponding to free groups
given in \cite[Proposition 3.7]{farsiseaton1}, note that the maps $e_{i_1, \ldots , i_l}$ and $I$ defined in
\cite[page 80]{ademleidaruan} are special cases of the construction above.  In
particular, let the free group $\F_k$ have generators $\gamma_1,
\ldots, \gamma_k$, and pick $\{ i_1, \ldots i_l \} \subseteq \{ 1,
2, \ldots , k \}$.  Define the homomorphism $\Phi_{i_1, \ldots ,
i_l} : \F_l \rightarrow \F_k$ by $\gamma_1 \mapsto \gamma_{i_1}$, and
then $e_{i_1, \ldots , i_l} = e_{\Phi_{i_1, \ldots , i_l}}$.
Similarly, let $\Phi_I : \Gamma \rightarrow \Gamma$ be the
isomorphism $\Phi_I(\gamma) = \gamma^{-1}$.  Then $e_{\Phi_I} = I$.

The homomorphism $e_\Phi$ induces a well-defined map
$|e_\Phi| : |\mathcal{G}^\Gamma | \rightarrow | \mathcal{G}^\Lambda |$ on orbit
spaces.  Letting $(1)$ denote the $\approx$-class of the trivial
homomorphism from $\Gamma$ to a unit in $\mathcal{G}$,
$\tilde{Q}_{(1)}$ is obviously diffeomorphic to $Q$. Recall that the
map $\pi : \tilde{Q}_\Gamma \rightarrow Q$ was defined in
\cite{farsiseaton1} by $\pi(\mathcal{G}\phi_x) = \mathcal{G}x$.
Letting $\iota : \{ 1 \} \rightarrow \Gamma$, we see that $\pi =
|e_\iota|$ up to the identification of $\tilde{Q}_{(1)}$ with $Q$.

In \cite{farsiseaton1} following Lemma 3.6, it was stated
that the map $\pi$ is not an embedding of the $\Gamma$-sectors of
$Q$ into $Q$.  By this, it was meant that the restriction of $\pi$
to a $\Gamma$-sector is not generally injective unless the local
groups of $Q$ are abelian. Using the definition \cite[Definition
2.3]{ademleidaruan} of embeddings of orbifold groupoids, the
restriction of this map to a $\Gamma$-sector is in fact an
embedding.

\begin{lemma}
\label{lem-piembedding}

Let $\mathcal{G}$ be an orbifold groupoid and $\Gamma$ a finitely
generated discrete group.  Let $(\phi) \in T_Q^\Gamma$ and $\iota : \{ 1 \}
\rightarrow \Gamma$, and then the restriction of the map $e_\iota$
to $(\phi)$ is an embedding of orbifold groupoids.

\end{lemma}

\begin{proof}

Let $x \in G_0$.  Then $e_\iota^{-1}(x) \cap (\phi)$ is precisely
the $G_x$-conjugacy class of $\phi_x$.  The number of such conjugacy
classes is given by the index of $C_{G_x}(\phi_x)$ in $G_x$, and so
for each orbifold chart $\{ V_x, G_x, \pi_x \}$ for $Q$ at $x$, we
have that $e_\iota^{-1}(V_x) \cap (\phi)$ is given by
\[
    \coprod\limits_{\psi_x \in (\phi_x)_{G_x}} V_x^{\langle \psi_x
    \rangle}
    =
    G_x/C_{G_x}(\phi_x) \times V_x^{\langle \psi_x \rangle},
\]
where $(\phi_x)_{G_x}$ denotes the $G_x$-conjugacy class of $\phi_x$.
Moreover, $\mathcal{G}^\Gamma_{e_\iota^{-1}(V_x)}$ is given by
\[
    \coprod\limits_{\psi_x \in (\phi_x)_{G_x}} C_{G_x}(\phi_x) \ltimes V_x^{\langle \psi_x
    \rangle}
    =
    G_x \ltimes  \left(G_x/C_{G_x}(\phi_x) \times V_x^{\langle \psi_x \rangle}\right).
\]
As $|e_\iota| = \pi$ is clearly proper, it follows that the
restriction of $e_\iota$ to each $\Gamma$-sector is an embedding.

\end{proof}

Hence, each $\Gamma$ sector of $Q$ is a suborbifold.  We
also have the following.

\begin{lemma}
\label{lem-ephisummand}

Let $\mathcal{G}$ be a groupoid, let $\Lambda$ and $\Gamma$ be finitely
generated discrete groups, and let $\Phi : \Lambda \rightarrow \Gamma$ be a
group homomorphism.  If $\Phi$ maps $\Lambda$ into $\Gamma$ as a
direct factor, then $e_\Phi$ is a finite union of embeddings of
orbifold groupoids, and the induced map $|e_\Phi|$ on orbit spaces
is surjective.  In particular, $e_\Phi$ is an embedding of each
$\Gamma$-sector into $\tilde{Q}_\Lambda$.  If $\Phi$ is an
isomorphism, then $e_\Phi$ is an isomorphism of orbifold groupoids.

\end{lemma}

By virtue of  \cite[Proposition 3.7]{farsiseaton1}, this result
generalizes \cite[Proposition 4.1]{ademleidaruan}.

\begin{proof}

Suppose $\Gamma = \Lambda \times \Gamma^\prime$ for some group
$\Gamma^\prime$ and $\Phi : \lambda \mapsto (\lambda, 1)$.  For each
$\phi_x \in \mathcal{S}_\mathcal{G}^\Lambda$, $e_{\Phi_0}(\phi_x
\times 1) = \phi_x$.  It follows that $e_{\Phi 0}$ is
surjective on objects and hence on orbit spaces.

The map $e^{-1} : \mathcal{G}^\Gamma \rightarrow
(\mathcal{G}^\Lambda)^{\Gamma^\prime}$ given in Theorem
\ref{thrm-Gamma1timesGamma2formula} is an isomorphism, and hence
$\mathcal{G}^\Gamma$ corresponds to the $\Gamma^\prime$-sectors of
$\mathcal{G}^\Lambda$. Moreover, the map $\pi :
(\mathcal{G}^\Lambda)^{\Gamma^\prime} \rightarrow
\mathcal{G}^\Lambda$ clearly satisfies $\pi \circ e^{-1} = e_\Phi$.  By
Lemma \ref{lem-piembedding}, then, the restriction of $e_\Phi$ to
each $\Gamma$-sector is an embedding.

If $e_\Phi$ is an isomorphism, then clearly $e_\Phi$ and
$e_{\Phi^{-1}}$ are inverse groupoid homomorphisms.

\end{proof}

\begin{lemma}
\label{lem-ephisurj}

Let $\mathcal{G}$ be a groupoid, let $\Lambda$ and $\Gamma$ be finitely
generated discrete groups, and let $\Phi : \Lambda \rightarrow \Gamma$ be a group
homomorphism.  If $\Phi$ is surjective, then $e_\Phi$ is an
embedding of $\mathcal{G}^\Gamma$ into $\mathcal{G}^\Lambda$ whose
image consists  of entire connected components.

\end{lemma}

\begin{proof}

Suppose $\Phi$ is surjective.  Then for each $\phi_x \in
\mathcal{S}_\mathcal{G}^\Gamma$, $\Imt \phi_x = \Imt e_{\Phi
0}(\phi_x)$ as subgroups of $G_x$.  Therefore, given a linear chart
$\{ V_x, G_x, \pi_x \}$ for $Q$ at $x$, $V_x^{\langle \phi_x
\rangle} = V_x^{\langle e_{\Phi 0}(\phi_x) \rangle}$, and
$C_{G_x}(\phi_x) = C_{G_x}(e_{\Phi 0}(\phi_x))$.  Moreover, if
$\phi_x, \phi_y^\prime \in \mathcal{S}_\mathcal{G}^\Gamma$ satisfy
$e_{\Phi 0}(\phi_x) = e_{\Phi 0}(\phi_y^\prime)$, then for each
$\lambda \in \Lambda$, we have $\phi_x(\Phi(\lambda)) =
\phi_y^\prime(\Phi(\lambda))$.  It clearly follows that $x = y$, and
as $\Phi$ is surjective, $\phi_x = \phi_y^\prime$.  Therefore,
$e_{\Phi 0}$ is injective.

By Equation \ref{eq-ePhi0local}, $e_{\Phi 0}$ is locally given by
the embedding of $V_x^{\langle \phi_x \rangle}$ into $V_x^{\langle
e_{\Phi 0}(\phi_x) \rangle}$.  As $V_x^{\langle \phi_x \rangle} =
V_x^{\langle e_{\Phi 0}(\phi_x) \rangle}$, this map is a submersion.
Moreover, we have that the restriction
$\mathcal{G}^\Gamma|_{V_x^{\langle \phi_x \rangle}}$ is given by
$C_{G_x}(\phi_x) \ltimes V_x^{\langle \phi_x \rangle}$, while the
restriction $\mathcal{G}^\Lambda|_{V_x^{\langle e_{\Phi 0}(\phi_x)
\rangle}}$ is given by $C_{G_x}(e_{\Phi 0}(\phi_x)) \ltimes
V_x^{\langle e_{\Phi 0}(\phi_x) \rangle} = C_{G_x}(\phi_x) \ltimes
V_x^{\langle \phi_x \rangle}$.  The maps $e_{\Phi 0}$ and $e_{\Phi
1}$ are simply the identity maps.  We see that $e_\Phi$ restricts in
each chart to an isomorphism of groupoids. As $e_{\Phi 0}$ is an
injective immersion, the induced map $|e_\Phi|$ on orbit spaces is
clearly proper.  It follows that $e_\Phi$ is an embedding of
orbifold groupoids.

To show that the image of $|e_\Phi|$ consists of entire connected
components, suppose $\phi_x \in \mathcal{S}_\mathcal{G}^\Gamma$ and
$\psi_y \in \mathcal{S}_\mathcal{G}^\Lambda$ such that
$\mathcal{G}(\phi_x\circ\Phi)$ and $\mathcal{G}\psi_y$ are elements
of the same connected component of $\tilde{Q}_\Lambda$.  By
\cite[page 9]{farsiseaton1}, it follows that $\phi_x\circ\Phi
\approx \psi_y$ as elements of $\mathcal{S}_\mathcal{G}^\Lambda$, so
that $\phi_x\circ\Phi$ and $\psi_y$ are connected by a sequence of
local equivalences in $\mathcal{S}_\mathcal{G}^\Lambda$.
Whenever $\psi_y \stackrel{loc}{\curvearrowright} \psi_{y^\prime}$,
there is a $g \in G_1$ such that $g[\xi_x^y\circ
\psi_{y^\prime}(\lambda)]g^{-1} = \psi_y(\lambda)$ for each $\lambda
\in \Lambda$.  Therefore, $\Imt \psi_y$ and $\Imt \psi_{y^\prime}$ are
isomorphic subgroups of $G_y$ and $G_{y^\prime}$, respectively.
Hence, there exists a $\phi_y : \Gamma \rightarrow G_y$ such that
$\phi_y \circ \Phi = \psi_y$ if and only if there exists a
$\phi_{y^\prime} : \Gamma \rightarrow G_{y^\prime}$ such that
$\phi_{y^\prime} \circ \Phi = \psi_{y^\prime}$.  As this is true for
each of the local equivalences connecting $\phi_x \circ \Phi$ to
$\psi_y$ (regardless of their direction), it follows that there is a
$\phi_y$ such that $\phi_y \circ \Phi = \psi_y$, and so $\psi_y$ is
in the image of $e_{\Phi 0}$.  It follows that $|e_{\Phi}|$ maps
connected components to connected components.

\end{proof}

It follows that when $\Phi$ is surjective, $e_\Phi$ embeds
$\tilde{Q}_\Gamma$ into $\tilde{Q}_\Lambda$ as a suborbifold
consisting entirely of $\Lambda$-sectors of $Q$.  In particular,
each $\Gamma$-sector of $Q$ is diffeomorphic to a $\Lambda$-sector
via $|e_\Phi|$.


\section{The $\Gamma$-Euler-Satake Characteristics}
\label{sec-GammaESC}

Throughout this section, we assume that $Q$ is a
closed orbifold presented by the orbifold groupoid $\mathcal{G}$.


\subsection{Definitions and Relations to Other Orbifold Euler Characteristics}
\label{subsec-defGammaESC}

We begin with the following.

\begin{definition}
\label{def-eulerchars}

Let $Q$ be closed and $\Gamma$ a finitely generated discrete group. Let
\[
\begin{array}{rcl}
    \chi_\Gamma^{ES}(Q)
    &=& \sum\limits_{(\phi) \in T_Q^\Gamma} \chi_{ES}\left(\tilde{Q}_{(\phi)}\right)
            \\\\
    &=& \chi_{ES} \left( \tilde{Q}_\Gamma \right)
\end{array}
\]
denote the {\bf $\Gamma$-Euler-Satake characteristic}, the
Euler-Satake characteristic of the space of $\Gamma$-sectors of $Q$
and
\[
\begin{array}{rcl}
    \chi_\Gamma^{top}(Q)
    &=& \sum\limits_{(\phi) \in T_Q^\Gamma} \chi_{top}\left(\tilde{Q}_{(\phi)}\right)
            \\\\
    &=& \chi_{top} \left( \tilde{Q}_\Gamma \right)
\end{array}
\]
the {\bf $\Gamma$-Euler characteristic}, the usual Euler
characteristic of the (underlying space of the) space of
$\Gamma$-sectors of $Q$.  Note that both sums are finite when
$\Gamma$ is finitely generated by \cite[Lemma 2.9]{farsiseaton1}.

\end{definition}

By virtue of the results in \cite{farsiseaton2}, these definitions
generalize those given in \cite{tamanoi1} and \cite{tamanoi2} for
global quotients.  In the case that $Q$ is oriented (and in
particular satisfies the codimension-2 requirement),
$\chi_\Gamma^{ES}(Q)$ is the evaluation of the $\Gamma$-Euler-Satake
class of $Q$ (see \cite[Section 2.3]{farsiseaton1}).

The following is a direct consequence of Theorem \ref{thrm-Gamma1timesGamma2formula}.

\begin{corollary}
\label{cor-GESCG1timesG2}

Let $Q$ be a closed orbifold and
$\Gamma_1$ and $\Gamma_2$ finitely generated discrete groups.  Then
\begin{equation}
\label{eq-eulercharinduction}
    \chi_{\Gamma_1\times\Gamma_2}^{ES}(Q)
    =
    \chi_{\Gamma_1}^{ES}(\tilde{Q}_{\Gamma_2}).
\end{equation}

\end{corollary}

We also have the following.

\begin{corollary}
\label{cor-generalproductstatement}

Let $Q_1$ and $Q_2$ be closed orbifolds and $\Gamma$ a finitely generated discrete group.
Then
\begin{equation}
\label{eq-GESCmultiplic}
    \chi_\Gamma^{ES}(Q_1 \times Q_2)
    =
    \chi_\Gamma^{ES}(Q_1)\chi_\Gamma^{ES}(Q_2).
\end{equation}

\end{corollary}

\begin{proof}

We note that
\begin{equation}
\label{eq-GESCmultiplic}
\begin{array}{rcl}
    \chi_\Gamma^{ES}(Q_1 \times Q_2)
    &=&     \chi_{ES}\left(\widetilde{(Q_1)}_\Gamma \times \widetilde{(Q_2)}_\Gamma\right)      \\\\
    &=&     \chi_{ES}\left(\widetilde{(Q_1)}_\Gamma\right)\chi_{ES}\left(\widetilde{(Q_2)}_\Gamma\right)
                    \\\\
    &=&     \chi_\Gamma^{ES}(Q_1)\chi_\Gamma^{ES}(Q_2),
\end{array}
\end{equation}
where the first equation follows from Proposition
\ref{prop-productformula} and the second from Theorem
\ref{thrm-ESCmultiplic}.

\end{proof}

Applying these results to the case of quotient orbifolds, we obtain the following
generalization of \cite[Proposition 2-1]{tamanoi1}.

\begin{corollary}
\label{cor-generalproductstatementquotients}

For $i = 1, 2$, let $G_i$ be a compact, connected Lie group and $M_i$
a smooth, compact manifold on which $G_i$ acts smoothly and locally freely.
Then for each finitely generated discrete group $\Gamma$,
\[
    \chi_\Gamma^{ES}((M_1 \times M_2) \rtimes
    (G_1 \times G_2))
    =
    \chi_\Gamma^{ES}(M_1 \rtimes G_1)
    \chi_\Gamma^{ES}(M_2 \rtimes G_2)
\]
and
\[
    \chi_\Gamma^{top}((M_1 \times M_2) \rtimes
    (G_1 \times G_2))
    =
    \chi_\Gamma^{top}(M_1 \rtimes G_1)
    \chi_\Gamma^{top}(M_2 \rtimes G_2).
\]
If $G$ is a compact, connected Lie group, $M$ is a smooth, compact
$G$-manifold, and $\Gamma_i$ is a finitely generated discrete group
for $i = 1, 2$, then
\[
    \chi_{\Gamma_1 \times \Gamma_2}^{ES}(M\rtimes G)
    =
    \sum\limits_{(\phi)\in\mbox{\scriptsize HOM}(\Gamma_1, G)/G}
    \chi_{\Gamma_2}^{ES} \left( M^{\langle \phi \rangle} \rtimes
    C_G(\phi) \right),
\]
and
\[
    \chi_{\Gamma_1 \times \Gamma_2}^{top}(M\rtimes G)
    =
    \sum\limits_{(\phi)\in\mbox{\scriptsize HOM}(\Gamma_1, G)/G}
    \chi_{\Gamma_2}^{top} \left( M^{\langle \phi \rangle} \rtimes
    C_G(\phi) \right).
\]

\end{corollary}

In the proof of \cite[Theorem 3.2]{seaton1}, it was
shown that $\chi_{ES}(\tilde{Q}) = \chi_{top}(Q)$, where $\tilde{Q}$
denotes the inertia orbifold (see \cite[Definition 2.49]{ademruan}).  Note
that this proof does not require the codimension-2 condition nor the
hypotheses that $Q$ is effective or orientable.  In
\cite[Corollary 3.8]{farsiseaton2}, it is shown that the inertia orbifold
corresponds to the $\Z$-sectors.  Hence, we have the relationship
\begin{equation}
\label{eq-eulercharstopinert}
    \chi_{top}(Q)
    =
    \chi_\Z^{ES}(Q)
    =
    \chi_{ES}(\tilde{Q}).
\end{equation}
More generally,
\begin{equation}
\label{eq-eulercharstopsatake}
    \chi_\Gamma^{top}(Q)
    =
    \chi_{\Gamma \times \Z}^{ES}(Q).
\end{equation}

If $Q$ is a manifold, then the $\Gamma$-sectors of $Q$ are simply
$Q$ so that $\chi_\Gamma^{ES}(Q) = \chi_{top}(Q)$ reduces to the
usual Euler characteristic for any $\Gamma$. In fact, the $\Gamma$-Euler-Satake
characteristic generalizes several Euler characteristics for
orbifolds.
\begin{itemize}
\item   When $\Gamma = \{ 1 \}$, $\tilde{Q}_\Gamma = Q$ so that
        $\chi_\Gamma^{top}(Q)$ is the usual Euler
        characteristic of the underlying topological space and $\chi_\Gamma^{ES}(Q)$ is the usual
        Euler-Satake characteristic.

\item   When $\Gamma = \Z^d$ and $Q$ is a global quotient, $\chi_\Gamma^{top}(Q)$ is the
        $d$th and $\chi_\Gamma^{ES}(Q)$ is the $d-1$st orbifold Euler characteristic defined in
        \cite{bryanfulman}.  This is seen by repeated application of Equations \ref{eq-eulercharinduction}
        and \ref{eq-eulercharstopinert}.  By Theorem
        \ref{thrm-Gamma1timesGamma2formula} and the fact that the
        $\Z$-sectors correspond to the inertia orbifold, this implies that the Euler characteristics
        of Bryan and Fulman are the Euler characteristics and Euler-Satake characteristics
        of the $d$th (respectively $d-1$st) inertia orbifold of $Q$.

\item   When $Q = M/G$ is a global quotient orbifold,
        $\chi_\Gamma^{top}(Q)$ and $\chi_\Gamma^{ES}(Q)$ are the
        generalized orbifold Euler characteristics defined in
        \cite{tamanoi1} and \cite{tamanoi2} (where they are denoted
        $\chi_\Gamma(M; G)$ and $\chi_\Gamma^{orb}(M; G)$,
        respectively).  This follows from \cite[Theorem 3.1]{farsiseaton2}.

\item   We have that $\chi_{\Z^2}^{ES}(Q) = \chi_{\Z}^{top}(Q)$ is the
        stringy orbifold Euler characteristic defined in
        \cite{dixon} for global quotients and generalized to general
        orbifolds in \cite{roan}.  Note that this Euler
        characteristic was shown to be the Euler characteristic of
        equivariant $K$-theory in \cite{segal} for quotients and for
        orbifold $K$-theory in \cite{ademruan} for quotient
        orbifolds; see also \cite{hirzebruchhoefer}.
        It follows that for each $d \geq 2$,
        $\chi_{\Z^d}^{ES}(Q) = \chi_{\Z^{d - 1}}^{top}(Q)$
        is the Euler characteristic of the orbifold
        $K$-theory of $\tilde{Q}_{\Z^{d - 2}}$.
\end{itemize}


\subsection{The Vanishing of the $\Gamma$-Euler-Satake Characteristics}
\label{subsec-vanishingGammaESC}

In \cite{farsiseaton1}, it was demonstrated that the Euler-Satake
characteristics of the individual $\Gamma$-sectors act as a complete
obstruction to the existence of a nonvanishing vector field on
closed codimension-2 orbifolds when $\Gamma$ covers the local groups
of $Q$.  Recall that this means that for each subgroup $H$ of an
isotropy group of $Q$, there is a surjective homomorphism from
$\Gamma$ onto $G$. Note that all of the $\Gamma$-Euler-Satake
characteristics may vanish even in the case that the
Euler-Satake characteristics of the individual sectors do not.

\begin{example}
\label{ex-totalvanishnotindivid}

Let $Q$ be the orbifold described as follows.  The underlying space of $Q$ is
homeomorphic to $\mathbb{T}^4$.  The singular strata of $Q$ consists of two connected components,
one diffeomorphic to $S^2$ and one diffeomorphic to the surface $\Sigma_2$ of genus $2$,
both modeled on $\R^2 \times (\R^2/(\Z/2\Z))$ where $\Z/2\Z$ acts as rotations.
Let $a$ generate $\Z/2\Z$ in a chart for the component diffeomorphic to $S^2$, and let
$b$ generate $\Z/2\Z$ in a chart for the component diffeomorphic to $\Sigma_2$;
clearly, the $\approx$-classes of $a$ and $b$, respectively, do not depend on the choice of chart.
Then it is easy to see that $\chi_{ES}(Q) = 0$.
Moreover, note that for any finitely generated discrete $\Gamma$, there is a
bijection between $\mbox{HOM}(\Gamma, \langle a \rangle)$ and
$\mbox{HOM}(\Gamma, \langle b \rangle)$, implying that there is a
$k$ such that
\[
\begin{array}{rcl}
    \chi_\Gamma^{ES}(Q)
        &=&     \chi_{ES}(Q) + k\chi_{ES}(S^2 \rtimes \langle a \rangle )
                + k\chi_{ES}(S^2 \rtimes \langle b \rangle )             \\\\
        &=&     0.
\end{array}
\]
However, for any $\Gamma$-sector corresponding to a homomorphism
with image $\langle a \rangle$ or $\langle b \rangle$, the Euler-Satake characteristic
of the corresponding $\Gamma$-sector is $1$ or $-1$, respectively.  Hence, it is possible
that each $\Gamma$-Euler-Satake characteristic vanishes even though the
Euler-Satake characteristic of the individual $\Gamma$-sectors do
not.

\end{example}

Here, we briefly note some consequences of the vanishing of the Euler-Satake
characteristic of the individual $\Gamma$-sectors of a closed orbifold $Q$.  First, we
have the following.

\begin{proposition}
\label{prop-nonewdiffeoclassesaftercoverslocalgroups}

Let $Q$ be a closed orbifold presented by the orbifold groupoid
$\mathcal{G}$ and suppose $\Gamma$ covers the local groups of $Q$.
Then for each finitely generated discrete group $\Lambda$ and each
$(\phi) \in T_Q^\Lambda$, there is a $(\psi) \in T_Q^\Gamma$ such
that $\tilde{Q}_{(\phi)}$ is orbifold-diffeomorphic to
$\tilde{Q}_{(\psi)}$.

\end{proposition}

\begin{proof}

Pick a representative $\phi_x$ of $(\phi)$, and as $\Gamma$ covers
the local groups of $Q$, there is a $\psi_x : \Gamma \rightarrow
G_x$ such that $\Imt\psi_x = \Imt\phi_x$.  Let $(\psi)$ denote the
$\approx$-class of $\psi_x$, and then it follows directly from the
definitions that there is a diffeomorphism between $(\phi)$ and
$(\psi)$.  It follows that $\tilde{Q}_{(\phi)}$ and
$\tilde{Q}_{(\psi)}$ are diffeomorphic.

\end{proof}

\begin{proposition}
\label{prop-esvanishifftopvanish}

Let $Q$ be a closed orbifold and $\Gamma$ a finitely generated
discrete group that covers the local groups of $Q$.  Suppose the
Euler-Satake characteristic of each $\Gamma$-sector is zero. Then
the topological Euler characteristic of (the underlying space of)
each $\Gamma$-sector is zero.

\end{proposition}

\begin{proof}

Suppose the Euler-Satake characteristic of each $\Gamma$-sector is
zero, and let $(\phi) \in T_Q^\Gamma$.  Then
\[
    \chi_{top}(\tilde{Q}_{(\phi)})
        =
    \chi_{ES}\left(\widetilde{\left(\tilde{Q}_{(\phi)} \right)}_\Z
    \right).
\]
Note that $\widetilde{\left(\tilde{Q}_{(\phi)} \right)}_\Z$ is a
disjoint union of $\Gamma \times \Z$-sectors of $Q$, and by
Proposition \ref{prop-nonewdiffeoclassesaftercoverslocalgroups},
each is diffeomorphic to a $\Gamma$-sector.  It follows that each
has Euler-Satake characteristic zero, so that
$\chi_{top}(\tilde{Q}_{(\phi)}) = 0$.

\end{proof}

The following is a consequence of Propositions
\ref{prop-nonewdiffeoclassesaftercoverslocalgroups} and
\ref{prop-esvanishifftopvanish}.

\begin{corollary}
\label{cor-equivtoallvanish}

Let $Q$ be a closed orbifold, and let $\Gamma$ be a finitely
generated discrete group that covers the local groups of $Q$.
Suppose $\Gamma$ admits a presentation with $m$ generators.  Then
the following are equivalent.
\newcounter{Acount}
\begin{list}{\roman{Acount}.}{\usecounter{Acount}}
\item   The Euler-Satake characteristic of each $\Gamma$-sector is
        $0$.
\item   For any finitely generated discrete group $\Gamma^\prime$
        that covers the local groups of $Q$, the
        Euler-Satake characteristic of each $\Gamma^\prime$-sector is
        zero.
\item   For each $d \geq m$, the Euler-Satake characteristic of each $\F_d$-sector is
        $0$.
\end{list}
Moreover if any one (hence all) of the above equivalent statements
is true, then the following hold.
\newcounter{Bcount}
\begin{list}{\arabic{Bcount}.}{\usecounter{Bcount}}
\item   For any finitely generated discrete group $\Lambda$, the
        Euler-Satake characteristic of each $\Lambda$-sector is zero.
\item   For any finitely generated discrete group $\Gamma^\prime$
        that covers the local groups of $Q$, the
        topological Euler characteristic of (the underlying space of)
        each $\Gamma^\prime$-sector is $0$.
\item   For each $d \geq m$, the topological Euler characteristic of (the underlying space of)
        each $\F_d$-sector is $0$.
\end{list}

\end{corollary}


\section{Formulas for Wreath Products}
\label{sec-wreathproducts}

In this section, we use the properties of wreath symmetric products and the $\Gamma$-Euler-Satake
characteristics demonstrated above to derive formulas for $\Z^m$-Euler-Satake characteristics of
orbifolds presented as wreath symmetric products by compact Lie groups.


\subsection{The Euler-Satake Characteristic of a Wreath Symmetric Product}
\label{subsec-ESCwreathprod}

Here, we derive the following formula for the Euler-Satake characteristic
of a wreath symmetric product $MG(\mathcal{S}_n)$ as defined in
Subsection \ref{subsec-actionswreathprodmanifolds}.

\begin{proposition}
\label{prop-ESCwreathprod}

Let $G$ be a compact, connected Lie group and $M$
a smooth, compact manifold on which $G$ acts smoothly and locally freely.
Then
\[
    \chi_{ES}(MG(\mathcal{S}_n))
    =
    \frac{1}{n!} \left(\chi_{ES}(Q)\right)^n.
\]

\end{proposition}

\begin{proof}

First, we show that the natural map $G^n x \mapsto
G(\mathcal{S}_n)x$ from the orbifold presented by $M^n \rtimes G^n$ to the
orbifold $MG(\mathcal{S}_n)$ is an orbifold cover of with $n!$
sheets.

Assume that $G$ acts effectively on $M$.  Let $p \in
MG(\mathcal{S}_n)$ be the orbit of a point $x \in M^n$. By
conjugating by elements of $G(\mathcal{S}_n)$, we can assume that
$x$ is of the form
\[
    x   =
    \left(\stackrel{r_1}{\overbrace{x_1, x_1, \dots, x_1}},\;
    \stackrel{r_2}{\overbrace{x_2, x_2, \dots x_2}},\;
    \ldots,\;
    \stackrel{r_N}{\overbrace{x_N, x_N, \dots, x_N}}
    \right)
\]
with each $x_j$ distinct, so that $\sum_{j=1}^N r_j = n$.

By Proposition \ref{prop-wrprodfixedptspecific}, the isotropy group
of $p$ with respect to the $G(\mathcal{S}_n)$-action is
\[
    (G({\mathcal S}_n ))_p \cong \prod\limits_{j= 1}^N  G_{x_j}({\mathcal S}_{r_j} ).
\]
Similarly, it is easy to see that the isotropy group of $p$ with respect to the $G^n$-action is
\[
    (G^n )_p = \prod_{j= 1}^N  (G_{x_j})^{r_j}.
\]
Note that $(G^n)_p$ is a normal subgroup of $G(\mathcal{S}_n)_p$,
and
\[
    (G(\mathcal{S}_n))_p/(G^n)_p
    \cong
    \prod\limits_{j=1}^N \mathcal{S}_{r_j} .
\]

For each $j$, choose a linear slice $U_{x_j}$ of the $G$-action on
$M$ centered at $x_j$, and let
\[
    W_p = \prod_{j=1}^N  (U_{x_j})^{r_j}.
\]
Then $W_p$ is a linear slice at $x$ for both the $G^n$ and
$G(\mathcal{S}_n)$-action on $M^n$.  Moreover,
\[
    W_p \times_{G(\mathcal{S}_n)_p} G(\mathcal{S}_n)
    \cong
    \left( W_p \times_{(G^n)_p} G^n \right) \times
    \mathcal{S}_n/ \left(\prod\limits_{j=1}^N \mathcal{S}_{r_j}
    \right),
\]
so that each point in $MG(\mathcal{S}_n)$ is the image of
$\left|\mathcal{S}_n/ \left(\prod\limits_{j=1}^N \mathcal{S}_{r_j}
\right)\right|$ points in $|M^n \rtimes G^n|$. In particular, if $x$
is nonsingular so that $G_x = \{ 1 \}$ and each of the $x_j$ are
distinct, then there are $|\mathcal{S}_n| = n!$ points in $|M^n
\rtimes G^n|$.   It follows that $M^n \rtimes G^n$ presents an
$n!$-sheeted cover of the orbifold presented by $MG(\mathcal{S}_n)$.

If $G$ does not act effectively on $M$, then let $K \trianglelefteq
G$ denotes the kernel of the action.  We apply the above argument
to the effective orbifolds presented by $G^n \rtimes K^n$ and
$M(G/K)(\mathcal{S}_n)$.  The kernel of the action
in in $M^n \rtimes G^n$ and $MG(\mathcal{S}_n)$ is $K^n$, and so the
result is the same.

Now, note that $M^n \rtimes G^n = (M\rtimes G)^n$
so that by Theorem \ref{thrm-ESCmultiplic}, $\chi_{ES}(M^n \rtimes G^n) = (\chi_{ES}(Q))^n$. Then the result follows from Lemma
\ref{lem-ESCmultipliconcovers}.

\end{proof}


\subsection{Generalizations of MacDonald's Formulas}
\label{subsec-macdonald}

In this section, we prove Theorem \ref{thrm-macdonaldformulas} which
generalizes MacDonald's formulas in \cite{macdonald} and
\cite[Theorem 5]{wang1} to the context of an orbifold given by
a quotient by a compact, connected Lie group acting smoothly and locally freely.
This formula will serve as the base case in the induction proof of
Theorem \ref{thrm-mainWPformula}.

We again assume that $G$ acts effectively for simplicity.  When this is not the case, we can always
apply these results to the associated effective orbifold as in the proof of Proposition \ref{prop-ESCwreathprod}. Then Equation \ref{eq-ESCnoneff} can be used to address the
the noneffective case.

\begin{theorem}
\label{thrm-macdonaldformulas}

Let $G$ be a compact, connected Lie group acting smoothly,
effectively, and almost freely on the closed manifold $M$.  Let $Q$ denote the
orbifold presented by $M \rtimes G$ and $MG({\mathcal S}_n)$ denote the wreath symmetric product.

\begin{enumerate}

\item
The cohomology $H^\ast(MG({\mathcal S}_n))$
satisfies the MacDonald formulas
\[
    H^\ast(MG({\mathcal S}_n)) \cong {\mathcal{SP}}^n [H^\ast(Q)]
\]
and
\[
    \sum_{n \geq 0} \mbox{dim}\: [H^\ast(MG({\mathcal S}_n))]q^n
    =
    \frac1{(1-q)^{\mbox{\scriptsize dim}\:[H^\ast(Q)]}}
\]
where $\mathcal{SP}$ denotes the symmetric product and
and $\mbox{dim}\: [W] = |\mbox{dim}\:[W^{ev}]\otimes \C|
- |\mbox{dim}\:[W^{od}]\otimes \C|$ denotes the Euler characteristic of the complex $W$.

\item
The delocalized equivariant cohomology $H_\Z^\ast(MG( {\mathcal S}_n))$ satisfies the MacDonald formula
\[
    \bigoplus_{n \geq 0}
    H_\Z^\ast(MG( {\mathcal S}_n) ) q^n
    \cong
    \mathcal{SP} \left[\bigoplus_{n \geq 0} H^\ast\left(\widetilde{(M/G)}_\Z\right)q^n\right].
\]
As well, we have the dimension formula
\[
    \sum_{n \geq 0} \mbox{dim} \: [H_\Z^\ast(MG({\mathcal S}_n))]q^n
    =
    \prod_{n\geq 1} \frac1{(1-q^n)^{\mbox{\scriptsize dim}\:(H_\Z^\ast(M \rtimes G))}}.
\]

\end{enumerate}
\end{theorem}

\begin{proof}

See \cite[Section 6]{wang1}, \cite{zhou}, and \cite[Section
4]{tamanoi1}.

To prove (1), by Equation \ref{eq-doublecrossprod} and the definition of
symmetric algebra, we have
\[
\begin{array}{rcl}
    H^\ast(MG( {\mathcal S}_n))
    &\cong& H^\ast(M^n)^{G^n\rtimes {\mathcal S}_n}
            \\\\
    &\cong&
    {\mathcal{SP}}^n [H^\ast(Q)].
\end{array}
\]
The dimension formula follows from the general fact that if  $V=
\bigoplus_p V^p$ is a general graded vector space,
\[
    \sum\limits_{n \geq 0} \mbox{dim}\:[\mathcal{SP}[V]] q^n
    =
    \frac{1}{(1-q)^{\mbox{\scriptsize dim} \: [V]}}.
\]

To prove (2), we compute the MacDonald formula for the delocalized
equivariant cohomology of $MG({\mathcal S}_n)$. We have that
\[
\begin{array}{rcl}
    \bigoplus_{n \geq 0} H_\Z^\ast(MG({\mathcal S}_n)) q^n
        &=&     \bigoplus\limits_{n \geq 0} q^n H^\ast\left(\widetilde{(MG({\mathcal S}_n)}_\Z\right)
                    \\\\
        &=&     \bigoplus\limits_{n \geq 0} q^n
                \bigoplus\limits_{(a) \in t_{M^n; G( \mathcal{S}_n)}} H^\ast\left((M^n)^{\langle a \rangle}
                \rtimes (C_{G({\mathcal S}_n)}(a))\right)
                    \\\\
        &\cong&
                \bigoplus\limits_{n \geq 0}  q^n \bigoplus\limits_{m_r(c)} \bigotimes\limits_{c, r}
                H^\ast \left((M^{\langle c\rangle})^{m_r(c)}\rtimes ((C_G(c)\langle a_{r,c} \rangle)^{m_r(c)}
                \rtimes {\mathcal S}_{m_r(c)})\right) .
\end{array}
\]
where  the sum over $m_r(c)$ is over all of the
$m_r(c)>0$ for $(c) \in t_{M; G}$ and $r = 1, \dots, n$ subject to
the constraint $\sum_{r,c} r m_{r}(c) = n$. By Equation
\ref{eq-doublecrossprod}, this is isomorphic to
\[
    \bigoplus\limits_{n \geq 0}  \bigoplus\limits_{m_r(c)}
    \bigotimes\limits_{c, r} q^{r m_r(c)} {\mathcal {SP}}^{m_r(c)}
    \left[H^\ast \left(M^{\langle c\rangle}\rtimes (C_G(c))\right)\right].
\]
Recall that $m_r = \sum_c m_r (c)$, and then this is isomorphic to
\[
    \bigoplus\limits_{n \geq 0}
    \bigoplus\limits_{m_r} {\mathcal {SP}}^{m_r}
    \left[ \bigoplus\limits_r H^\ast\left(\widetilde{(M\rtimes G)}_\Z\right) q^r\right],
\]
where the sum over $m_r$ is the sum over all of the $m_r>0$ such
that $\sum_{r} r m_{r} =n$. Therefore, if we assume $r\geq 1$ and
$m_r\geq 1$ with no other constraint, this is isomorphic to
\[
    \bigoplus\limits_{m_r} {\mathcal {SP}}^{  m_r}
    \left[\bigoplus\limits_r H^\ast\left(\widetilde{(M\rtimes G)}_\Z\right) q^r\right].
\]

By the definition of the symmetric product, this last expression is
equal to
\[
    \mathcal{SP}\left[\bigoplus_r H^\ast\left(\widetilde{(M\rtimes G)}_\Z\right)
    q^r\right].
\]
Note that this formula only holds for $r, m_r = 1,2, \dots$,
and not for $r,m_r = 1,2, \dots, n$.

For the dimension formula, we apply \cite[Lemma 3.1]{wangzhou} to
\[
    \mathcal{SP} \left[\bigoplus\limits_{r \in \mathbb{N}}  H^\ast\left(\widetilde{(M\rtimes G)}_\Z\right)
    q^r\right]
    \cong
    \bigotimes\limits_{r \in \mathbb{N}} \mathcal{SP} \left[ H^\ast\left(\widetilde{(M\rtimes G)}_\Z\right) q^r\right].
\]

\end{proof}

The $K$-theory of the crossed product algebra $C(M) \rtimes G$ is isomorphic to the equivariant
$K$-theory of $M$ (see e.g. \cite{blackadar}); i.e.
\[
\begin{array}{rcl}
    K^i(C^\ast(M\rtimes G))
        &\cong&     K^i(C(M) \rtimes G)         \\\\
        &=&         K^i_G(C(M))                 \\\\
        &=&         K_i^G(M)
\end{array}
\]
for $i = 1, 2$.  For an orbifold  groupoid $\mathcal{G}$, recall from Section \ref{subsec-defGammaESC}
that $\chi_{\Z^2}^{ES}(Q)$ coincides with the stringy orbifold Euler characteristic and the
Euler characteristic of the orbifold $K$-theory of $\mathcal{G}$; i.e.
\[
    \chi_{\Z^2}^{ES}(Q)
    =
    \mbox{dim}\;K^0(C^*(\mathcal{G}))\otimes \C - \mbox{dim}\; K^1(C^*(\mathcal{G}))\otimes \C
\]
(see \cite[Definition 3.8]{ademleidaruan}).  Moreover, in the case of a quotient orbifold
$\mathcal{G} = M \rtimes G$, the (complexified) $K$-theory is isomorphic to the delocalized cohomology.
Hence, we have the following.

\begin{corollary}

With $M$ and $G$ as in Theorem \ref{thrm-macdonaldformulas},
\[
    \sum_{n \geq 0} \chi_{\Z^2}^{ES}(MG(\mathcal{S}_n)) q^n
    =
    \sum_{n \geq 0} \mbox{dim}\: [K^\ast(MG( {\mathcal S}_n))]q^n
    =
    \prod_{n\geq 1} \frac1{(1-q^n)^{\mbox{\scriptsize dim}\:(K^\ast(M \rtimes G))}}.
\]

\end{corollary}


\subsection{Wreath Symmetric Products and Product Formulas}
\label{subsec-wreathproductformulas}

In this section, we prove a product formula for the $\Z^n$-Euler
characteristics that extends \cite[Theorem A]{tamanoi1} to the case
of a quotient by a compact, connected group.  For simplicity, we
abreviate $\chi_{\Z^m}^{top}$ as $\chi_{(m)}$ and
$\chi_{\Z^m}^{ES}$ as $\chi_{(m)}^{ES}$.

The following follows \cite[Lemma 4-1]{tamanoi1}.

\begin{lemma}
\label{lem-ECwithtrivaction}

Let $L$ be a compact Lie group, $K$ a closed connected subgroup, and
suppose there is an $a \in L$ such that $\langle a \rangle \cap K =
\langle a^r \rangle$ and $a \in C_L(K)$.  Let $M$ be a $K$-manifold on which $K$ acts
smoothly, effectively, and locally freely, and suppose
$a$ acts trivially on $M$.  Then for each $m \geq 0$,
\[
    \chi_{(m)} (M\rtimes K\langle a\rangle )
    =
    r^m \chi_{(m)} (M\rtimes K).
\]

\end{lemma}

\begin{proof}

Note that $K\langle a\rangle = \bigcup_{i=1}^r S a^i$.  The case $m
= 0$ is obvious, so suppose $m \geq 1$.

By definition,
\[
    \chi_{(m)} (M \rtimes K\langle a \rangle)
    =
    \sum_{(\phi) \in t_{M;K\langle a \rangle}^{\Z^m}}
    \chi_{top}\left( M^{\langle \phi \rangle} \rtimes C_{K\langle a \rangle}
    (\phi)\right).
\]
It is straightforward that the centralizer of $\phi$ in
$K\langle a \rangle$ is equal to $C_K (\phi)C_{\langle a \rangle}
(\phi)$.  Since $C_{\langle a \rangle} (\phi)$ is a subgroup of
$\langle a \rangle$, it acts trivially on $M$.  It follows that for
each $(\phi) \in t_{M;K\langle a \rangle}^{\Z^m}$,
\[
    \chi_{top}\left( M^{\langle \phi \rangle} \rtimes C_{K\langle a \rangle} (\phi)\right)
    =
    \chi_{top}\left( M^{\langle \phi \rangle} \rtimes C_K (\phi)\right).
\]

The image of $\phi$ is a set of $m$ commuting elements of $S\langle
a \rangle$.  Note that for $g, h \in K$, $g a^i$ and $h a^j$ commute if
and only if $g$ and $h$ commute.  Moreover, for each commuting $m$-tuple
$g_1, \dots, g_m$ of elements of $K$, there are $r^m$ corresponding
$m$-tuples $g_1a^{i_1}, \dots, g_ma^{i_m}$ of commuting elements of
$K\langle a \rangle$.  It follows that
\[
    \chi_{(m)} (M\rtimes K\langle a \rangle)
    =
    r^m \sum_{(\phi)\in t_{M;K\langle a \rangle}^{\Z^m}}
    \chi_{top} \left( M^{\langle \phi \rangle} \rtimes C_K (\phi) \right),
\]
finishing the proof.

\end{proof}

\begin{theorem}
\label{thrm-mainWPformula}

Let $G$ be a compact, connected Lie group acting smoothly,
effectively, and almost freely on the closed manifold $M$, and let $Q$
be the orbifold presented by $M\rtimes G$.  Then
\begin{equation}
\label{eq-MacdonaldExp}
    \sum\limits_{n=0}^{\infty} \chi_{ES}(MG(\mathcal{S}_n)) q^n
    =
    \exp \left[ q\chi_{ES}(Q) \right].
\end{equation}
For $m \geq 0$,
\begin{equation}
\label{eq-genfunctionwreath}
    \sum\limits_{n=0}^{\infty} \chi_{(m)}(MG(\mathcal{S}_n)) q^n
    =
    \prod\limits_{j_1, \dots, j_{m} \geq 1}
    \left[ \left( 1 - q^{j_1 j_2\dots j_m} \right)^{j_2 j_3^2 \dots j_{m-1}^{m-2} j_m^{m-1}}
    \right]^{- \chi_{(m)}(Q)}.
\end{equation}
Noting that $J_{r,m} = \sum_{j_1 \dots j_m =r} \ j_2 j_3^2 \dots
j_{m-1}^{m-2} j_m^{m-1} $ is the number of subgroups of index $r$ in
$\Z^m$, Equation \ref{eq-genfunctionwreath} can be rewritten as
\begin{equation}
\label{eq-genfunctionwreathsimplified}
    \sum\limits_{n=0}^{\infty} \chi_{(m)}(MG(\mathcal{S}_n))q^n
    =
    \prod\limits_{r \geq 1}
    \left[\left( 1-q^r \right)^{J_{r,m}} \right]^{-\chi_{(m)}(Q)}.
\end{equation}

\end{theorem}

\begin{proof}

Equation \ref{eq-MacdonaldExp} follows by applying Proposition \ref{prop-ESCwreathprod}.
The proof of Equations \ref{eq-genfunctionwreath} and \ref{eq-genfunctionwreathsimplified}
is by induction on $m$.  As $\chi_{(0)}(Q) = \chi_{top}(Q)$
and $\chi_{(1)}(Q) = \chi_{top}(\tilde{Q}_\Z)$, the
cases $m = 0$ and $m = 1$ of Equation \ref{eq-genfunctionwreath}
corresponds to
\[
    \sum\limits_{n=0}^{\infty} \chi_{top}(MG(\mathcal{S}_n))q^n
    =
    \left( 1-q \right)^{-\chi_{top}(Q)}
\]
and
\[
    \sum\limits_{n=0}^{\infty} \chi_{top}\left( \widetilde{MG(\mathcal{S}_n)}_\Z \right)q^n
    =
    \prod\limits_{r \geq 1} \left( 1-q^n \right)^{-\chi_{top}(\tilde{Q}_\Z)},
\]
both following directly from Theorem \ref{thrm-macdonaldformulas}.

For $m > 1$, we have by applying Proposition
\ref{prop-homsplitsurjectiveandcentralizer} that
\[
\begin{array}{rcl}
    \sum\limits_{n\geq 0} q^n \chi_{(m)}(MG({\mathcal S}_n))
    &=&
    \sum_{n\geq 0} q^n \chi_{(m-1)}\left(
    \widetilde{(MG(\mathcal{S}_n))_\Z}
    \right)             \\\\
    &=&
    \sum\limits_{n\geq 0} q^n
    \sum\limits_{(({\mathbf g }, s )) \in t_{M^n; G(\mathcal{S}_n)}} \chi_{(m-1)}
    (M^{\langle {\mathbf g }, s \rangle}
    \rtimes C_{G({\mathcal S}_n)} ({\mathbf g }, s )).
                        \\\\
    &=&
    \sum\limits_{n \geq 0} q^n
    \chi_{(m-1)} \left(
    \coprod\limits_{m_r(c)}
    \prod\limits_r
    \prod\limits_{(c)}
    \left( (M^{\langle c \rangle})^{m_r(c)}
    \rtimes  ([C_G(c)\langle a_{r,c}\rangle]^{m_r(c)} \rtimes {\mathcal S}_{m_r(c)} ) \right)\right),
\end{array}
\]
where the union and products are over all $m_r(c) > 0$ for $(c) \in t_{M;G}$,
$r\geq 1$ subject to the constraint that $\sum_{r,c} rm_r(c) = n$,
and $(c) \in t_{M;G}$. Applying Corollary \ref{cor-generalproductstatementquotients}, we have
\[
\begin{array}{rcl}
    \sum\limits_{n\geq 0} q^n \chi_{(m)}(MG({\mathcal S}_n))
    &=&
    \sum\limits_{n \geq 0} q^n \sum\limits_{m_r(c)}
    \prod\limits_r
    \prod\limits_{(c)}
    \chi_{(m-1)}\left( (M^{\langle c \rangle})^{m_r(c)}
    \rtimes  ([C_G(c)\langle a_{r,c}\rangle]^{m_r(c)} \rtimes {\mathcal S}_{m_r(c)} ) \right),
                        \\\\
    &=&
    \sum\limits_{n\geq 0} q^n
    \sum\limits_{m_r(c)} \prod\limits_r
    \prod\limits_{(c)}
    \chi_{(m-1)}\left( \mathcal{SP}^{m_r (c)}(M^{\langle c\rangle}) \rtimes
    (C_G(c)\langle a_{r,c}\rangle) \right),
\end{array}
\]
where the symmetric product of a crossed product groupoid $X = K
\rtimes L$  is defined in a similar way to the symmetric product of
a vector space; that is, ${\mathcal {SP}}^{m}[X]= K^m \rtimes ({L^m
\rtimes \mathcal S}_m)$. In other words, the symmetric product is
the wreath product.

As in the proof of Theorem \ref{thrm-macdonaldformulas}, the above
formula implies that if we again set $m_r = \sum_c m_r(c)$, we have
\[
\begin{array}{rcl}
    \sum\limits_{n\geq 0} q^n \chi_{(m)}(MG({\mathcal S}_n))
    &=&
    \sum\limits_{n\geq 0}
    q^n \chi_{(m-1)}
    \left( \widetilde{(MG({\mathcal S}_n))}_\Z \right)
                    \\\\
    &=&
    \sum\limits_{n\geq 0}
    q^n  \sum\limits_{m_r} \prod\limits_r
    \chi_{(m-1)}\left( {\mathcal{SP}}^{m_r}[{\mathbb{IMG}_r}]
    \right),
\end{array}
\]
where now the sum is over all $m_r > 0$ and $r \geq 1$ such that
$\sum_r r m_r = n$, and $\mathcal{SP}^{m_r}= \bigoplus_{(c)}
\mathcal{SP}^{m_r(c)}$.  Note that the groupoid $\mathbb{IMG}_r$ is
defined in Equation \ref{eq-wreathinertiagroupoid}. Hence, as in the
proof of Theorem \ref{thrm-macdonaldformulas} (2), this is equal to
\[
    \sum\limits_{m_r} \prod\limits_r \chi_{(m-1)} \left(
    \mathcal{SP}^{m_r} [\mathbb{IMG}_r q ] \right)
\]
with no constraints on $m_r \geq 0$.  As $\chi_{(m-1)}$ is a
multiplicative character, this is equal to
\[
    \prod\limits_r \sum\limits_{m_r} \chi_{(m-1)} \left(
    \mathcal{SP}^{m_r} [\mathbb{IMG}_r q ] \right).
\]

Now, by the inductive hypothesis,
\begin{equation}
\label{eq-WPformulainductivehyp}
\begin{array}{c}
    \sum\limits_{n\geq 0} q^n \chi_{(m)}(MG({\mathcal S}_n))
    \\
    =
    \prod\limits_r \left[ \prod\limits_{r_1,\dots, r_{m-1} \geq 1}
    \left( 1-(q^r)^{{r_1,\dots, r_{m-1} }} \right)^{r_2 r_3^2\dots r_{m-1}^{m-2}}\right]^{-\chi_{(m-1)}
    ({\mathbb{IMG}}_r)}.
\end{array}
\end{equation}
By Lemma \ref{lem-ECwithtrivaction}, Theorem \ref{thrm-centeralizerwreathprod}, and
Proposition \ref{prop-wrprodfixedptgeneral},
\[
    \chi_{(m-1)}(\mathbb{IMG}_r) = r^{m-1}\chi_{(m-1)}\left(
    \widetilde{(M \rtimes G)}_\Z
    \right).
\]
Combining this with Equation \ref{eq-WPformulainductivehyp}, we have
\[
\begin{array}{rcl}
    \sum\limits_{n\geq 0} q^n \chi_{(m)}(MG({\mathcal S}_n))
    &=&
    \left[ \prod\limits_{ r, r_1,\dots, r_{m-1} \geq 1}
    \left( 1-(q^r)^{{r_1,\dots, r_{m-1} }} \right)^{r_2 r_3^2\dots
    r_{m-1}^{m-2}}\right]^{-  r^{m-1}\chi^{(m-1)}\left(
    \widetilde{(MG(\mathcal{S}_n))}\right)}
            \\\\
    &=&
    \left[ \prod\limits_{r, r_1,\dots, r_{m-1} \geq 1}
    \left( 1-q^{{r_1,\dots, r_{m-1}r }} \right)^{r_2 r_3^2\dots
    r_{m-1}^{m-2}r^{m-1}}\right]^{-\chi_{(m)}(Q)},
\end{array}
\]
completing the proof.

\end{proof}

An application of Equation \ref{eq-eulercharstopsatake} allows us to express
Theorem \ref{thrm-mainWPformula} in terms of the $\Z^m$-Euler-Satake characteristics as
follows.

\begin{corollary}
\label{cor-mainWPformulaESC}

Let $G$ be a compact, connected Lie group acting smoothly,
effectively, and almost freely on the closed manifold $M$, and let $Q$
be the orbifold presented by $M\rtimes G$.  Then
\[
    \sum\limits_{n=0}^{\infty} \chi_{(m)}^{ES}(MG(\mathcal{S}_n)) q^n
    =
    \left\{
    \begin{array}{ll}
    \exp \left[ q\chi_{ES}(Q) \right],  &   m = 0,
    \\\\
    \prod\limits_{r \geq 1}
    \left[\left( 1-q^r \right)^{J_{r,m+1}} \right]^{-\chi_{(m)}^{ES}(Q)},
    &   m > 0.
    \end{array} \right.
\]

\end{corollary}


\section{Generalized Hodge Numbers}
\label{sec-hodgenumbers}

In this section, we will generalize to arbitrary orbifolds a product
formula that Wang and Zhou proved for orbifold Hodge numbers of
global quotient orbifolds in \cite{wangzhou}; see Theorem \ref{thrm-hodgeproductformula}.
Note that for special choices of the parameters $x$ and $y$ in Theorem \ref{thrm-hodgeproductformula},
we obtain interesting geometric invariants; see the note after Definition \ref{def-hodgepolysstandard}.

To derive our product formula for the shifted Hodge numbers, we will
use the product formula for the case of the Hodge Poincar\'{e}
numbers and make the necessary modifications that account for the
degree shifting.  We assume that $G$ is a compact complex connected
Lie group acting effectively, holomorphically, and almost freely on
the compact complex manifold $M$, and $Q$ is the complex orbifold
presented by $M \rtimes G$. See Subsection \ref{subsec-cohomology}
for the definition of the orbifold Dolbeault cohomology of $M
\rtimes G$. The following definition follows \cite[page
157]{wangzhou}.

\begin{definition}
\label{def-hodgepolysstandard}

Let the compact, complex, connected Lie group $G$ act effectively,
locally freely, and holomorphically on the complex manifold $M$. We
define the \emph{standard orbifold Hodge polynomials} $H(M \rtimes
G; x,y)$ and $h(M\rtimes G; x,y)$ by
\[
    H(M\rtimes G; x,y)
    =
    \sum\limits_{s,t \geq 0} H^{s,t}(M\rtimes G)x^s y^t
\]
and
\[
    h(M \rtimes G; x,y)
    =
    \sum\limits_{s,t \geq 0}
    \left|\mbox{dim}_\C [H^{s,t}(M\rtimes G)]\right|x^s y^t,
\]
where $H^{s,t}(M \rtimes G)$ are the standard Dolbeault cohomology
groups of the complex orbifold $M\rtimes G$; see Subsection
\ref{subsec-cohomology}.  We define the \emph{standard
delocalized orbifold Hodge polynomials} $H_D(M\rtimes G; x,y)$ and
$h_D(M\rtimes G; x,y)$ by
\[
    H_D(M\rtimes G; x,y)
    =
    \sum\limits_{s,t \geq 0}
    H^{s,t}\left( \tilde{Q}_\Z \right) x^s y^t,
\]
and
\[
    h_D(M\rtimes G; x,y)
    =
    \sum\limits_{s,t \geq 0}
    \left|\mbox{dim}_\C [H^{s,t}\left( \tilde{Q}_\Z \right)\right| x^s y^t .
\]

\end{definition}

For example, it is well known that the Hirzebruch genus of $M\rtimes G$ is given by
$h(M\rtimes G, y, -1)$, the (topological) Euler characteristic by
$h(M\rtimes G, -1, -1)$, and the signature by $h(M\rtimes G, 1, -1)$.

In the remaining definitions, we will need to consider the grading
shift.  In this context, the degree shifting number of \cite[Section 4.2]{ademleidaruan}
corresponds to that of \cite{zaslow}, which we now recall.

Fix $c \in G$, and let $M^{\langle c \rangle}_1, \dots M^{\langle c
\rangle}_{N_c}$ be the connected components of $M^{\langle c
\rangle}$. We let $F^c_j$ denote the shift number associated to each
$M^{\langle c\rangle}_j$, $j=1, \dots, N_c$, as follows. As $c$
fixes $M^{\langle c\rangle}_j$, the action of $c$ on the tangent
space of each point in  $M^{\langle c\rangle}_j$ can be represented
by a diagonal unitary matrix
\[
    \mbox{Diag}
    \left(e^{2\pi i \theta_1^j}, \dots, e^{2\pi i \theta_{d_j}^j}\right)
\]
with $d_j$ the complex dimension of $M^{\langle c \rangle}_j$, $0 <
\theta_i^j \leq 1$ for each $i = 1, \dots, d_j$ and $j = 1, \dots,
N_c$. Define the shift number
\[
    F_j^c
    =
    \sum\limits_{k=1}^{d_j} \theta^j_k \in \Q.
\]
Then for $0 \leq p, q \leq \mbox{dim}_\C \: Q$, the Chen-Ruan
Dolbeault cohomology groups are given by
\[
    H_{CR}^{p,q}(M \rtimes G)
    =
    \bigoplus\limits_{(c) \in t_{M;G}} \bigoplus\limits_{j = 1}^{N_c}
    H^{p - F_j^c,q - F^c_j}
    \left(M^{\langle c \rangle}_j \rtimes C_G(c)\right).
\]

For the remainder of this section, we restrict to the case of $G$
acting on $M$ such that all shift numbers are integers.  The
following definitions follow \cite[page 156]{wangzhou}.

\begin{definition}
\label{def-hodgepolysshifted}

Let the compact, complex, connected Lie group $G$ act locally freely
and holomorphically on the compact, complex manifold $M$ and suppose
each of the shift numbers $F_j^c$ are integers.  Define the
\emph{shifted delocalized orbifold Hodge polynomials}
$h_{CR}(M\rtimes G; x,y)$ by
\[
    h_{CR}(M\rtimes G; x,y)
    =
    \sum\limits_{s,t \geq 0}
    \left|\mbox{dim}\: [H_{CR}^{s,t} (M\rtimes G) \right| x^s y^t .
\]
Also define the \emph{shifted delocalized orbifold Hodge numbers}
$\mathcal{H}^{p,q}(M\rtimes G)$ of $M\rtimes G$ for $0 \leq p, q \leq
\mbox{dim}_\C\: Q$  by
\[
    {\mathcal H}^{p,q}(M \rtimes G)
    =
    \sum\limits_{(c) \in t_{M;G}} \;\;
    \sum\limits_{j = 1}^{N_c}
    h_j^{p-F_j^c, q - F^c_j}(c),
\]
where
\[
\begin{array}{rcl}
    h^{p,q}_j(c)
    &=&
    \mbox{dim}_\C \left[H^{p,q}\left(M_j^{\langle c\rangle} \rtimes C_G(c)\right)\right]
                \\\\
    &=&
    (-1)^{p+q} \mbox{dim}_\C \left[H^{p,q}\left(M_j^{\langle c\rangle} \rtimes C_G(c)\right)\right],
\end{array}
\]
for $j=1, \dots, N_g$, $p, q \geq 0$.

\end{definition}

Abstractly, given a bigraded vector space $V = \bigoplus_{s, t \geq
0} V^{s, t}$, we can define a new grading on $V$ by using any
integer shifts.  Clearly, this shift will be reflected at the
dimension level, thus giving rise to formulas similar to those in
\cite[page 5]{zhou}.  In particular, we will use the formal notation
$V\{\{K\}\}$ for $\bigoplus_{s, t \geq 0} V^{s - K, t - K}$.

To derive our product formula, we will use the following, which
follows from general results on graded vector spaces (see \cite[page
163]{wangzhou}).  In particular, the proof follows that of Theorem
\ref{thrm-macdonaldformulas} (2).

\begin{proposition}
\label{prop-hodgedelocalized}

Let $N$ be a compact, complex manifold and let $K$ be a compact
complex, connected Lie group acting effectively, holomorphically,
and locally freely on $N$ with integer shifts.

\begin{enumerate}

\item   The standard delocalized orbifold Hodge polynomials $H_D$ satisfy
\[
    \bigoplus\limits_{n \geq 0} H_D(N^n \rtimes (K^n\rtimes {\mathcal S}_n); x,y) q^n
    =
    \mathcal{SP}\left[\bigoplus\limits_{r>0}  H_D(N\rtimes K; x,y) q^r\right],
\]
where $\mathcal{SP}$ again denotes the symmetric product algebra.

\item   The standard delocalized orbifold Hodge polynomials $h_D$ satisfy
\[
    \sum\limits_{n \geq 0} h_D (N^n \rtimes (K^n\rtimes {\mathcal S}_n) ; x,y)q^n
    =
    \prod\limits_{n, s,t}\frac{1}{ (1- x^s y^t q^n)^{h_D^{s,t}(N \rtimes K) }},
\]
where
\[
\begin{array}{rcl}
    h_D^{s,t}(N \rtimes K)
    &=&
    \mbox{dim}_\C \left(H_D^{s,t}(N \rtimes K)\right)
                \\\\
    &=&
    (-1)^{s + t}|\mbox{dim}_\C (H_D^{s,t}(N \rtimes K) |,
\end{array}
\]
for $s,t \geq 0$, and where
\[
    h_D(N\rtimes K) =  \sum\limits_{s,t}h_D^{s,t}(N\rtimes K).
\]

\end{enumerate}

\end{proposition}

The following is the main result of this section.
For the case of $G$ finite, see \cite[Theorem 3.1]{wangzhou}.

\begin{theorem}
\label{thrm-hodgeproductformula}

Let $M$ be a compact, complex manifold on which the compact, complex, connected Lie group $G$ acts effectively, holomorphically, and locally freely with integer shifts.
The shifted delocalized orbifold Hodge polynomials satisfy
\[
    \sum_{n \geq 0} h_{CR}(MG({\mathcal S}_n); -x, -y) q^n
    =
    \prod\limits_{n=1}^{\infty} \;\prod\limits_{s,t\geq 0} \;
    \frac{1}{ (1- x^s y^t q^n (xy)^{(r-1)d/2})^{(-1)^{s+t} h_{CR}^{s,t}(M\rtimes G) }}.
\]

\end{theorem}

\begin{proof}

The result follows from Proposition \ref{prop-hodgedelocalized} adjusted for the appropriate
degree shift as in \cite[page 5]{zhou}, together with Lemma \ref{lem-hodgeshiftedwreath} below.

\end{proof}

\begin{lemma}
\label{lem-hodgeshiftedwreath}

Let $M$ be a compact, complex manifold on
which the compact, complex, connected Lie group $G$ acts effectively,
holomorphically, and almost freely with integer shifts. Then the
degree shifts of the wreath product orbifold $MG(\mathcal{S}_n)$
are given by
\[
    F_{\rho}
    =
    \sum\limits_{r = 1}^n
    \sum\limits_{(c) \in t_{M;G}} \sum\limits_{j = 1}^{N_c}
    m_{r,(c)}(j) \left(F^c_j+ d \frac{r-1}{2}\right),
\]
where $\rho = \{ m_r(c) \}_{r \geq 1, (c) \in t_{M;G}}$,
$m_r(c) = \sum_{j=1}^{N_c} m_{r,(c)}(j)$, and $\sum_{r,(c)} rm_r(c) = n$.

\end{lemma}

\begin{proof}

By using the results of Subsection \ref{subsec-wreathprodcompact},
we will compute the shifts for the wreath product orbifold
$NK({\mathcal S}_n)$ presented by $N^n \rtimes (K^n\rtimes {\mathcal S}_n)$ in
terms of the shifts of the orbifold presented by $N\rtimes K$.  This yields the local
situation on $MG(\mathcal{S}_n)$.  By Proposition \ref{prop-cycledecompwrprod}, it is
enough to consider an
element $W$ of $NK({\mathcal S}_n)$ of the form $W = ( (g,1, \dots,
1), (12\dots n))$, with  $g \in K$. By direct calculation, a fixed
point in $N^n$ of the action of $W$ is of type $(x,\dots, x)$, where
$x\in N^{\langle g\rangle}$. Since the calculation can be made
locally, we assume that we have local coordinates $(z_1, \dots,
z_d)$ near $x\in N^{\langle g\rangle}$ with the action of $K$ given
by
\[
    g(z_1, \dots, z_d)
    =
    ( e^{2\pi i \theta_1}z_1, \dots,
    e^{2\pi i \theta_d} z_d),
\]
with $\theta_j \in [0,1)$ and $j=1, \ldots d$.  Hence $g\in K $ is
locally given by
\[
    \mbox{Diag}\left( e^{2\pi i \theta_1}, \dots, e^{2\pi i \theta_d}\right),
\]
where we can assume that $\theta_{r+1}, \ldots , \theta_d =0$.
Then near $(x,\dots, x)\in N^n$, $x \in N^{\langle g\rangle}$, $W$ is
given by a block--diagonal
matrix with blocks of the form
\[
\left[  \begin{array}{cccccc}
    0 &0 &0 &\dots&0&e^{2\pi i \theta_j}\\
    1 &0 &0 &\dots&0&0\\
    0 &1 &0 &\dots&0&0\\
    0 &0 &1 &\dots&0&0\\
    \dots &\dots &\dots &\dots&\dots&\dots\\
    0 &0 &0&\dots &1&0\\
\end{array}
\right].
\]
A straightforward calculation shows that the
characteristic polynomial of the above matrix is
\[
    \lambda^n - e^{2\pi i \theta_j},
\]
and so its eigenvalues are given by
\[
    e^{2\pi i \frac{(\theta_j+k)}n}
\]
for $k = 0,\dots, n-1$.

It follows that the shift $F^W(x,\dots, x)$ at $(x, \dots, x)$ for
the component of $(N^n)^W$ containing $(x, \dots, x)$, $x \in N$, is
given by
\[
F^W(x,\dots, x) = \sum\limits_{j=1}^r \sum\limits_{k=0}^{n-1}\frac{{(\theta_j+k)}}n + (d-r)
\sum\limits_{k=0}^{n-1} \frac{k}n,
\]
where the first sum represents the terms arising from the
eigenvalues $\not= 1$, and the second from the eigenvalues $=1$.
Using elementary algebra, the previous equation becomes
\[
\begin{array}{rcl}
    F^W(x,\dots, x)
    &=&
    \sum\limits_{j=1}^r \frac{\theta_j}n + d\, \frac{n-1}{2}
                    \\\\
    &=&
    F_j^c + d \frac{n-1}{2},
\end{array}
\]
where $F_j^c$ is the shift associated to the component $M^{\langle
c\rangle}_j$ of the fixed set $N^{\langle g\rangle}$ for the action
of $C_K(c)$.

Now consider a general conjugacy class containing an element of type
$\rho = \{  m_r(c) \}_{r \geq 1, (c) \in t_{N;K}}$, where
$\sum_{r,(c)} r m_r(c) =n$. The description of the fixed set
$(M^n)^{\langle a\rangle}$ given in Proposition
\ref{prop-wrprodfixedptspecific} implies that the $\rho$-component
can be written  as
\begin{equation}
\label{eq-hodgefixedpointset}
    (M^n)^{\langle a\rangle}_\rho
    =
    \prod\limits_{r,(c)} \prod_{j =1}^{N_c}
    {\mathcal {SP}}^{m_{r,(c)}(j)} [M^{\langle c\rangle}_j \rtimes C_G(c)],
\end{equation}
where $\mathcal{SP}$ denotes the symmetric product and
$(m_{r,(c)}(1), \dots, m_{r,(c)}(N_c))$ satisfies $\sum_{j=1}^{N_c}
m_{r,(c)}(j) = m_r(c)$. Then the shift $F_{\rho}$ corresponding to the connected
component described in Equation
\ref{eq-hodgefixedpointset} is given by
\[
    F_{\rho}
    =
    \sum\limits_{r,c}
    \sum\limits_{j = 1}^{N_c} m_{r,(c)}(j)
    (F^c_j + d \frac{r-1}2).
\]

\end{proof}


\bibliographystyle{amsplain}

\end{document}